\documentclass[10pt]{amsart}
\usepackage{epsfig, psfrag, amscd}

\newcommand{\R}{{\mathbb R}}  
  
\newcommand{\I}{{\mathbb I}}    
\newcommand{\Z}{{\mathbb Z}}
\newcommand{\Op}{{\mathcal O}}  
\newcommand{\Fo}{{\text{FSet}^{op}}}
\newcommand{\Set}{{\text{FSet}}}
\newcommand{\F}{{E}}

\newcommand{\la}{\langle [}
\newcommand{\ra}{] \rangle}  
\newcommand{\n}{{\mathbf{n}}}

\newcommand{\und}{\underline}
\newcommand{\wt}{\widetilde} 
\newcommand{\mC}{{\mathcal{C}}}
\newcommand{\Emb}{{\rm Emb}}
\newcommand{\Imm}{{\rm Imm}}

\newcommand{\mB}{{\mathcal B}}

\newcommand{\one}{{\bf{1}}}
\newcommand{\K}{{\mathcal K}}
\newcommand{\A}{{\mathcal A}}
\newcommand{\D}{{\mathcal D}}
\newcommand{\G}{{\mathcal G}}
\newcommand{\Top}{{\bf{\rm Top}}}

\newcommand{\U}{{\mathcal{U}}}
\def\holim#1#2{\mathchoice
	{\mathop{\displaystyle\mathop{\mbox{\rm holim}}_{\longleftarrow}}
	_{#1}{#2}}
	{\mathop{\smash{\displaystyle\mathop{\mbox{\rm holim}}
	_{\longleftarrow}}}
	_{#1}{#2}\vphantom{\displaystyle\mathop{\mbox{\rm holim}}
	_{\longleftarrow}}}
	{}{}}
\newcommand{\hofib}{{\rm{hofib}}\;}

\newcommand{\mF}{{\mathcal{F}}}
\newcommand{\tot}{\widetilde{\rm Tot}}

\newcommand{\refT}[1]{Theorem~\ref{T:#1}}                          
\newcommand{\refC}[1]{Corollary~\ref{C:#1}}                          
\newcommand{\refP}[1]{Proposition~\ref{P:#1}}                          
\newcommand{\refD}[1]{Definition~\ref{D:#1}}                          
\newcommand{\refL}[1]{Lemma~\ref{L:#1}}                          
\newcommand{\refE}[1]{Equation~\ref{E:#1}} 
\newcommand{\refF}[1]{Figure~\ref{F:#1}}         
\newcommand{\refX}[1]{Example~\ref{X:#1}}    
\newcommand{\refS}[1]{Section~\ref{S:#1}}

\theoremstyle{plain}                          
\newtheorem{theorem}{Theorem}[section]                          
\newtheorem{proposition}[theorem]{Proposition}                          
\newtheorem{lemma}[theorem]{Lemma}                          
\newtheorem{corollary}[theorem]{Corollary}                          
 
\newtheorem{mydiagram}[theorem]{Figure}
\newtheorem{example}[theorem]{Example}

\theoremstyle{definition}                          
\newtheorem{definition}[theorem]{Definition}  

\theoremstyle{remark}                          
\newtheorem*{remark}{Remark}                          
\newtheorem*{notation}{Notation}

\begin{document}

\title{Operads and knot spaces}
\author{Dev P. Sinha}
\address{Department of Mathematics, University of Oregon, Eugene, OR 97403}
\thanks{The  author is partially supported by NSF grant DMS-0405922.}

\maketitle

\section{Introduction}

Let $\F_m$ denote the space of embeddings of the interval $\I = [-1,1]$
in the cube $\I^m$ with endpoints and tangent vectors at those endpoints
fixed on opposite faces of the cube, equipped with a homotopy
through immersions to the unknot -- see \refD{embimm}.    
By \refP{null}, $\F_m$ is homotopy
equivalent to $\Emb(\I, \I^m) \times \Omega \Imm(\I, \I^m)$.  In \cite{McSm01}, McClure
and Smith define a cosimplicial object $\Op^\bullet$ associated to an operad
with multiplication $\Op$, whose homotopy invariant
totalization we denote $\tot(\Op^\bullet)$ --
see \refD{hoch} and \refD{tot} below.   Let $\K_m$ denote
the $m$th Kontsevich operad, introduced in \cite{Kont99},
whose entries are compactified configuration spaces and which 
is weakly equivalent to the little $m$-disks operad \cite{Sinh03} -- see 
\refD{cnrm} and \refT{kop} below.  

\begin{theorem}\label{T:main}
The totalization of the Kontsevich operad $\tot(\K^\bullet_m)$
is homotopy equivalent to the inverse limit of the Taylor tower approximations 
for $\F_m$ in the calculus of embeddings.  Moreover, $\tot^n(\K_m^\bullet)$
is the $n$th degree approximation.
\end{theorem}

Building on work of Goodwillie, Klein and Weiss \cite{Weis99, GoWe99, Good03, GoKl03},
and Volic \cite{Voli04.1, Voli04.2}, we have the following.

\begin{corollary} \label{C:main}
For $m>3$, $\F_m$ is weakly equivalent to $\tot(\K^\bullet_m)$.  
For $m=3$, all real-valued finite-type invariants of framed knots factor 
through a map from $\F_m$ to $\tot(\K^\bullet_m)$.
\end{corollary}

Applying the homology spectral sequence of a cosimplicial space, we have the following.

\begin{corollary}\label{C:HHss}
For $m>3$, there is a spectral sequence with $E^2$ page given by the 
Hochschild cohomology of the degree $m-1$ Poisson operad 
and which converges to the homology of $\F_m$.
\end{corollary}

These results resolve conjectures of Kontsevich from his address 
at the AMS Mathematical Challenges Conference 
in the summer of 2000  \cite{Kont00}.  Kontsevich's insights were motivated by
novel combinatorial work of Tourtchine.  In \cite{Tour03} Tourtchine gives
an algebraic description of the $E^1$-term of Vassiliev's homology
spectral sequence closely related to of our \refC{HHss}.  
Our results are at the level of spaces and show  that the disagreement 
which Tourtchine found between 
Hochschild cohomology of the Poisson operad  and 
Vassiliev's $E^2$-term is accounted for  by the fact that
$\F_m$ is not the classical knot space but is instead the
space of knots with trivialization through immersions. 

Our results bring together some recent 
developments in algebraic topology and its
application to fields such as deformation theory and knot theory.  In
\cite{Sinh02} we presented models for spaces of knots, including a
cosimplical model which is analogous to the cosimplicial model for loop
spaces.  We build on these results in proving \refT{main}.
In \cite{McSm01} McClure and Smith resolved the integral Deligne
conjecture, showing that the totalization of an operad with multiplication 
has a two-cubes action both in 
the setting of chain complexes and that of spaces.  
We apply their results to establish the following.

\begin{theorem}\label{T:2cube}
For any $m$, there is a little two-cubes action on $\tot(\K^\bullet_m)$.  For $m>3$,
$\F_m$ is a two-fold loop space.
\end{theorem}

We conjecture that this two-cubes action on $\tot(\K^\bullet_m)$ is compatible 
with a two-cubes action on the space of framed knots which
has been recently defined by Budney \cite{Budn03}, who 
goes on two show that long knots in dimension three are free over the two-cubes
action.  In future work we plan to investigate 
analogues of this freeness result in 
higher dimensions.  A first step will be to construct operations compatible
with this the two-cubes structure in the homology spectral sequence for 
an operad with multiplication, as McClure and Smith currently plan to do.
On  the $E^1$-term such operations will presumably coincide with Tourtchine's bracket,
defined combinatorially in \cite{Tour03}, 
 but through their space-level construction would also be compatible with
differentials and extend to further terms.

Some of the technical results developed in this paper may 
be of independent interest.  
We fully develop the operad structure, with multiplication, on the simplicial
compactification of configurations in Euclidean spaces.   An operad
structure on the canonical (Axelrod-Singer) compactification is known
\cite{GeJo94, Mark99} but does not yield an operad with multiplication.
Instead there is a map from Stasheff's
$A_\infty$ operad.  Our approach to the operad structure on the 
simplicial variant  blends geometry and combinatorics, 
revealing an operad structure on the 
standard simplicial model for the two-sphere.

\subsection{Acknowledgements}

We thank J. McClure and J. Smith for their interest in this project, answers
to questions, and
especially for writing Section~15 of \cite{McSm02}.  We thank M. Markl and
J. Stasheff for comments on early versions of this work, and M. Kontsevich 
for helpful conversations.

\tableofcontents

\section{Background material}

Our main results are stated in terms of operads with multiplication 
and their associated cosimplicial spaces.  As a chance to 
set the choice of definitions and notation which will be most convenient, 
and as an opportunity to place all standard material together,
we review this material here.  For a more complete survey
we highly recommend \cite{McSm04}.  In particular, 
Section~3 of \cite{McSm04} introduces cosimplicial spaces,
and Section~6 briefly introduces operads.  The paper  \cite{MSS02}
gives a more comprehensive 
introduction to operads.
A reader familiar with these
constructions may wish to
skip this section and refer back for clarification as needed.

\subsection{Cosimplicial spaces} \label{S:css}

\begin{definition}
Let $\Delta$ denote the category
with one object for every non-negative integer and
where the morphisms from $k$ to $\ell$ are the order-preserving maps
from $[k] = \{0, \cdots, k\}$ to $[\ell] = \{0, \cdots, \ell\}$, ordered 
in the standard way. 
A cosimplicial object in a category
${\mathcal{C}}$ is a (covariant) functor from $\Delta$ to ${\mathcal{C}}$.
A simplicial object is a contravariant functor from $\Delta$
to ${\mathcal{C}}$.
\end{definition}

Cosimplicial objects are denoted $X^\bullet$, where $X^k$ is the image
of $[k]$ under the functor, also known as the $k$th entry.
Simplicial objects are denoted $X_\bullet$.
Central in the theory is the standard cosimplicial space $\Delta^\bullet$, whose $k$th
entry is $\Delta^k$, with vertices labelled  by $[k]$, 
and which sends a morphism $[k] \to [\ell]$ to the linear map which extends
this map on vertices.

Every order-preserving map $[k] \to [\ell]$ can be factored through
elementary maps $d^i$, which are an isomorphism but for one element -
namely $i$ - not in their image, and elementary maps $s^i$, which 
are an isomorphism but for having $i$ and $i+1$ in $[k]$ 
both mapping to $i \in [\ell]$.   The corresponding maps between
entries of simplicial and cosimplicial objects, called (co)face and (co)degeneracy
maps, are often taken as a basis for their definition.
The definitions are arranged so that the simplices of an ordered simplicial
complex form a simplicial set.  Indeed, a simplicial set or simplicial
space $X_\bullet$ determines a space called its realization and 
denoted $|X_\bullet|$, defined as the
quotient space of the union of $X_i \times \Delta^i$ over all $i$ by
the relations $d_j x \times \beta \sim x \times d^j \beta$ and
$s_j x \times \beta \sim x \times s^j \beta$ for all $x \in X^i$, 
$\beta \in \Delta^i$.

Cosimplicial spaces naturally arise when studying mapping spaces.
The totalization of a cosimplicial space ${\rm Tot}X^\bullet$ 
is  the space of natural transformations from  $\Delta^\bullet$ to $X^\bullet$,
which is first used tautologically to study mapping spaces as follows.

\begin{definition} \label{D:xtoy}
For any $X \in \mC$, a category whose categorical product
$\odot$ is symmetric monoidal, taking the product of $X$ with itself gives
rise to a functor $X^{-} : \Fo \to \mC$ which sends $S$ to $\bigodot_{s \in S} X$, 
where $\Set$ is the category of finite sets.
By composing a simplicial set $Y_\bullet : \Delta \to \Fo$
with this functor, we obtain a cosimplicial object $X^{Y_\bullet}$.
\end{definition} 
 
\begin{proposition}
If $Y_\bullet$ is a simplicial set and $X$ is a space, ${\rm Tot}(X^{Y_\bullet})$ 
is homeomorphic to the space of maps from $|Y_\bullet|$ to $X$.  
\end{proposition}

For based $X$ and $Y_\bullet$ we may replace $X^{Y_n}$ by the 
subspace consisting of based maps from $Y_n$ to $X$.
The resulting cosimplicial space, which we denote $X_\star^{Y_\bullet}$, has
totalization homeomorphic to the space of based maps from
$|Y_\bullet|$ to $X$.
Another interesting example along these lines
is the Hochschild simplicial vector space $A^{S^1_\bullet}$,
whose associated chain complex computes Hochschild homology
of a commutative algebra $A$. 

Cosimplicial spaces are intimately connected with homotopy limits
(in fact, homotopy limits are defined in terms of cosimplicial
spaces in \cite{BoKa72}).
The nerve of a category $C$ is the simplicial set $C_\bullet$,
with $C_i$ being the collection of $i$ composable morphisms
and structure maps defined through composing such maps
or inserting identity maps (see for example Chapter~14 of \cite{Hirs03}).  
Denote the realization of the nerve of $C$ by $BC$, also called
the classifying space.  
Recall that if $c$ is an object of $C$, the category $C
\downarrow c$  has objects which are maps with target 
$c$ and morphisms given by morphisms in $C$ which commute with these 
structure maps.  The classifying space  $B(C \downarrow c)$ is contractible 
because $C \downarrow c$ has a final object, namely $c$  mapping to itself by the 
identity morphism.  A morphism $g$ from $c$ to $d$ induces a map  
from $C \downarrow c$ to $C \downarrow d$, so that $B(C\downarrow -)$ 
is itself a functor from $C$ to spaces. 

\begin{definition}
The homotopy limit of a functor $\F$ from a small category $C$ to the 
category of spaces is ${\rm Nat}(B(C \downarrow -), \F)$,
the space of natural transformations from $B(C \downarrow -)$ to $\F$.
\end{definition}

For $X^\bullet$ fibrant in the standard model structure on cosimplicial spaces, 
that is those which satisfy the matching condition 10.4.6 of  \cite{BoKa72}, 
 Theorem 11.4.4 of \cite{BoKa72} states 
 ${\rm Tot}(X^\bullet) \simeq \holim{\Delta}{X^\bullet}$, an equivalence needed for
 many applications.    For cosimplicial spaces which do not necessarily
satisfy the matching condition, we use an alternate definition of totalization
for which this equivalence is a tautology. 

\begin{definition}\label{D:tot}
\begin{itemize}
\item Let $\wt{\Delta^\bullet}$ be the cosimplicial space whose $[k]$th entry
is $B(\Delta \downarrow [k])$ and whose structure maps are the standard induced maps.
\item For a cosimplicial space $X^\bullet$ let $\tot X^\bullet$, called
the homotopy invariant totalization, denote the space
of natural transformations from $\wt{\Delta^\bullet}$ to $X^\bullet$.  
\item Let $\tot^k X^\bullet$
denote the space of natural transformations from the $k$th coskeleton
of $\wt{\Delta^\bullet}$ to $X^\bullet$.  
\item Let $\Delta_k$ denote the full subcategory of $\Delta$ whose objects
are $[i]$ for $i \leq k$.  Let $i_k :\Delta_k \to \Delta$ be the inclusion functor.
\end{itemize}
\end{definition}

In Section~15 of \cite{McSm02} the notations $\wt{\Delta^\bullet}$
and $\tot$ are used for any cofibrant replacement for $\Delta^\bullet$
and the corresponding totalization in the model structure on 
cosimplicial spaces where all objects are fibrant (in the usual
model structure from \cite{BoKa72}, all objects are cofibrant).
We choose one model of cofibrant replacement for definiteness.  

\begin{definition}
Define $\tot(X^\bullet)$ to be the homotopy limit of $X^\bullet$ and 
$\tot^k(X^\bullet) \cong \holim{}{(X^\bullet \circ i_k)}$.  
\end{definition}

The cosimplicial category $\Delta$
is also intimately related to the category of subsets of a finite set.

\begin{definition}\label{D:P0}
Let $P(k)$ be the category of all subsets of $[k] = \{0, 
\cdots, k\}$ where morphisms are defined by inclusion.  Let $P_0(k)$ be
the full subcategory of non-empty subsets. 
\end{definition}

The connection of this category to the simplicial world is evident in
the identification of $BP_0(k)$ with the barycentric subdivision of
a $k$-simplex.  We use this identification to
define maps $\wt{\Delta}^k \to \Delta^k$, and thus
${\rm Tot} X^\bullet \to \tot X^\bullet$ for any $X^\bullet$,
through the identification of $\wt{\Delta}^k$ and $\Delta^k$ with
$B(\Delta \downarrow [k])$ and $BP_0(k)$ respectively.
Namely, take the map induced on classifying spaces by the functor
which sends some $[n] \overset{f}{\to} [k]$ in $(\Delta \downarrow [k])$
to the image of $f$, as a subset of $[k]$.
There is also a translation between
cosimplicial diagrams and those indexed by $P_0(k)$, which we will use
in \refS{main}.

\begin{definition}\label{D:ck} 
Let $c_k \colon P_0(k) \to \Delta$ be the functor which sends a subset 
$S$ to the object in $\Delta$ with the same cardinality, and which
sends an inclusion $S \subseteq S'$ to the composite $[i]
\cong S \subset S' \cong [j]$, where $[i]$ and $[j]$ are
isomorphic to $S$ and $S'$ respectively as ordered sets.  
\end{definition}

The following lemma is a consequence of Theorem~6.4 of \cite{Sinh02}, 
along with the
observation that $\tot^k(X^\bullet) = \holim{}{(i_k \circ X^\bullet)}$.

\begin{lemma}\label{L:coshol}
For $X^\bullet$ a cosimplicial space,
$\holim{}{ \left(X^\bullet \circ c_k \right)}$ is weakly equivalent to $\tot^k X^\bullet$.
\end{lemma}

For a cosimplicial space there are spectral sequences for the
homotopy groups \cite{BoKa72} and homology groups \cite{Bous86, Ship96} of its 
(homotopy invariant) totalization, which we will apply in \refS{cons}.  
The homotopy spectral sequence is straightforward, with convergence 
immediate from its definition through the tower 
$$\cdots \leftarrow {\rm Tot}^i X^\bullet  \leftarrow {\rm Tot}^{i+1} X^\bullet \leftarrow \cdots,$$
whose homotopy inverse limit is ${\rm Tot} X^\bullet$ and which is a tower
of fibrations if $X^\bullet$ is fibrant.   Unraveling
the definitions we have the following.

\begin{proposition}\label{P:csshot}
For a fibrant cosimplicial space $X^\bullet$
there is a spectral sequence converging to $\pi_*({\rm Tot} X^\bullet)$
with $E^1_{-p,q} = \bigcap \; {\text{ker}} \; {s^k}_* \subseteq 
\pi_q(X^p). $ The $d_1$ differential is the restriction to 
this kernel of the map $\Sigma_{i=0}^{p+1} (-1)^i {d^i}_* \colon 
\pi_q(X^{p-1}) \to \pi_q(X^p).$ 
\end{proposition}

The homology spectral
sequence is more subtle in its convergence.
It generalizes the Eilenberg-Moore spectral sequence.
One of the precise statements  as to the convergence
of this spectral sequence arising from \cite{Bous86}
is as follows.  

\begin{theorem}\label{T:csshom}
For a fibrant 
cosimplicial space $X^\bullet$ there is a spectral sequence with
$E^1_{-p,q} = \bigcap \; {\text{ker}} \; {s^k}_* \subseteq 
H_q(X^p). $ The $d_1$ differential is the restriction to 
this kernel of the map 
$$\Sigma_{i=0}^{p+1} (-1)^i {d^i}_* \colon 
H_q(X^{p-1}) \to H_q(X^p).$$   
This spectral sequence
converges to $H_*({\rm Tot} X^\bullet)$ if  $X^k$ is simply
connected for all $k$ and $E^1_{-p,q} = 0$
when $q \geq cp$ for some $c>1$.  

Alternately, one
may arrive at the same spectral sequence from $E^2$ forward
with $E^1_{-p,q} = H_q(X^p)$ and $d_1$ defined as before,
but not restricted to the kernel of the codegeneracies.
\end{theorem}

This theorem follows immediately from 
Theorem 3.2 of \cite{Bous86} and the Universal
Coefficient Theorem, as both of Bousfield's
conditions, namely that 
$E^1_{-p,q} = 0$ if $p>q$ and that only finitely
many $E^1_{-p,q}$ with $q-p = n$ are non-zero for any given $n$,
follow from the vanishing with $q \geq cp$ for some $c>1$.   

These spectral sequences 
apply unchanged to the homotopy invariant totalization, in which case
the fibrancy condition can be dropped.  If $X^\bullet$
is a cosimpicial space and $\underline{X^\bullet}$ is a fibrant replacement 
(as given by Proposition~8.1.3 and Theorem~15.3.4 in \cite{Hirs03}) then 
\begin{equation*}
 \tot(X^\bullet) = {\rm Maps}\left(\wt{\Delta^\bullet}, X^\bullet \right)
\simeq {\rm Maps} \left( \wt{\Delta^\bullet}, \underline{X^\bullet} \right)
\simeq {\rm Maps} \left( \Delta^\bullet, \underline{X^\bullet} \right)
= {\rm Tot}(\underline{X^\bullet}).
\end{equation*}
Because homotopy and homology of the entries and structure maps of
$\underline{X^\bullet}$ agree with those of $X^\bullet$, the identifications of the
$E^1$-terms of the associated spectral sequences are unchanged.

\subsection{Operads}

We define non-$\Sigma$ operads in terms of a 
well-known \cite{Boar71, MSS02, Sinh02, Sinh03} category of rooted trees.

\begin{definition}
\begin{itemize}
\item A rooted, planar tree (or rp-tree) is an isotopy class of
finite connected acyclic graph with a distinguished vertex called the
root, embedded in the upper half plane with the root at the origin
so that the vertical coordinate in the plane is a monotone function
which increases on each edge as the distance from the root increases. 
Univalent vertices of an rp-tree (not counting the root, if it is 
univalent) are called leaves.

\item  Each edge of the tree is oriented by the  direction of the
root path, which is the unique shortest path to the root.  The vertex of
an edge which is further from the root is called its inital vertex, and
the vertex closer to the root is called its terminal vertex.  We say
that one vertex or edge lies over another if the latter is in the
root path of the former.  A non-root edge is called redundant if its 
initial vertex is bivalent.

\item Given an rp-tree $T$ and a set of edges $E$ the 
contraction of $T$ by $E$ is the tree $T'$ obtained by, for 
each edge $e \in E$, identifying its initial vertex with its
terminal vertex (altering the embedding of the tree
in a neighborhood of $e$) and removing $e$ from the set of
edges. 

\item Let $\Upsilon$ denote the category of rp-trees, where there is
a  morphism from $T$ to $T'$, denoted either $f_{T, T'}$ or $c_E$,
if $T'$ is the
contraction of $T$ along the set of  non-leaf edges $E$.  

\item Both the  collection of leaves in an rp-tree and the collections of 
edges with a given  terminal vertex  are ordered, 
using the clockwise orientation of the plane.

\item A sub-tree of an rp-tree is a connected sub-graph.  A sub-tree
is an rp-tree through a linear isotopy which translates $v$
to the origin.  
\end{itemize}
\end{definition}




See \refF{basicmor} for some examples of objects and morphisms in $\Upsilon$.
Let  $\Upsilon_n$ denote the full subcategory
of trees with $n$ leaves.  Note that $\Upsilon_n$ differs from $\Phi_n$ of
\cite{Sinh02}, which is canonically isomorphic to 
the full subcategory of rp-trees without bivalent vertices.  Each
$\Upsilon_n$ has a terminal object, namely the unique tree with one vertex,
called the $n$th corolla $\gamma_n$ as in \cite{MSS02}.  We allow for the tree
$\gamma_0$ which has no leaves, only a root vertex, and is
the only element of $\Upsilon_0$.
For a vertex $v$ let $|v|$ denote the number of edges for which
$v$ is terminal, usually called the arity of $v$.

\begin{definition}\label{D:op}
A non-$\Sigma$ operad is a functor $\Op$ from $\Upsilon$ to a symmetric
monoidal category $(\mC, \odot)$ which satisfies the following axioms.

\begin{enumerate}
\item $\Op(T) = \odot_{v \in T} \Op(\gamma_{|v|})$.  \label{1}
\item $\Op(\gamma_1) = \one_{\mC} = \Op(\gamma_0)$. \label{2}
\item If $e$ is a redundant edge and $v$ is its terminal vertex
then $\Op(c_{\{e\}})$ is the identity map on
$\odot_{v' \neq v} F(\gamma_{v'})$ tensored  with the isomorphism
$(\one_{\mC} \odot -)$ under the  decomposition of axiom (\ref{1}).  \label{3}
\item If $S$ is a subtree of $T$ and if $f_{S, S'}$ and $f_{T, T'}$
contract the same set of edges, then under the
decomposition of (\ref{1}), $F(f_{T,T'}) = F(f_{S, S'}) \odot id$. \label{4}
\end{enumerate}
\end{definition}

We sketch the equivalence of this definition with two standard ones.
By axiom~(\ref{1}), the values of $\Op$ are determined by 
its values on the corollas $\Op(\gamma_n)$, which corresponds
to $\Op(n)$ in the usual operad terminology of \cite{May72}.  
Axioms~(\ref{2}) and~(\ref{3})
correspond to the unit condition. By axiom~(\ref{4}),
the values of $\Op$ on morphisms may be computed by
composing morphisms on sub-trees, so we may identify some 
subset of basic morphisms through which all morphisms
factor.  In \refF{basicmor} we illustrate some basic morphisms in $\Upsilon$.
The first  corresponds to what are known as $\circ_i$ operations.
The second corresponds May's
operad structure maps from Definition~1.1 of \cite{May72}.  
That $\Op$ is a functor implies 
the commutativity of diagrams involving these basic morphisms.
Another basic class extending these is 
that of all morphisms $T \to \gamma_n$ where $\gamma_n$ is a 
corolla and $T$ is any tree.

\begin{center}
\begin{mydiagram}\label{F:basicmor}

\begin{center}

\psfrag{C1}{}
\psfrag{C2}{}
$$\includegraphics[width=12cm]{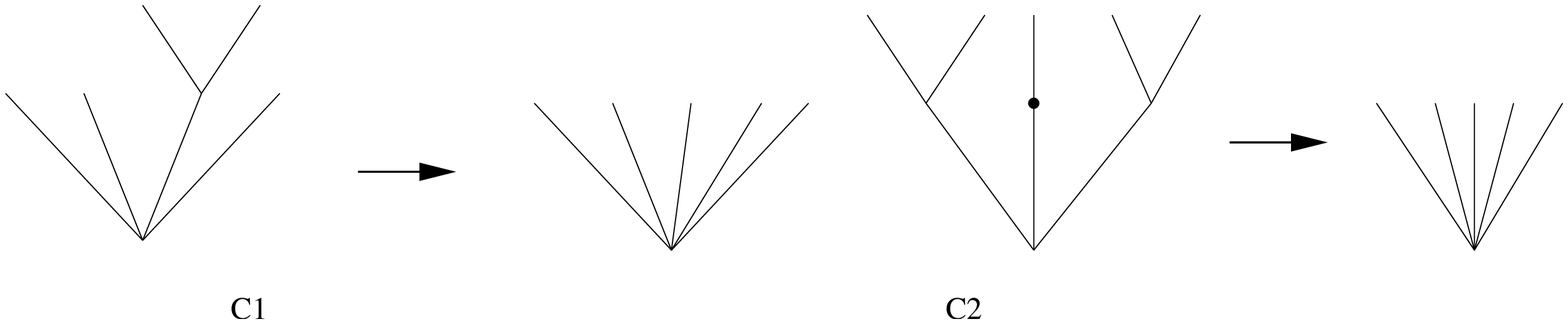}$$
Two morphisms in $\Upsilon$ which give rise to standard operad structure maps.\\
The first corresponds to a $\circ_i$ operation, the second to one of May's structure maps.
\end{center}
\end{mydiagram}

\end{center}

\begin{example}
\begin{itemize}
\item The associative operad $\A$, defined in any symmetric monoidal
category, has $\A(T) = \one_\mC$, and $\A(T \to T') = id$ for all morphisms 
in $\Upsilon$.

\item Let $\Phi$ denote the full subcategory of rp-trees with no redundant
edges (called the category of reduced trees in \cite{MSS02})
and let $P: \Upsilon \to \Phi$ denote the functor which contracts all of the
redundant edges of an rp-tree. 
The operad of planar trees, $\it{Tree}_n$ from
Definition~1.41 of \cite{MSS02},
is the operad in the category of
sets which sends $T$ to the set  of all
$T' \in \Phi$ which map to $P(T)$.
It sends a contraction of edges of $T$ to the collection of
contractions on the corresponding edges for trees over $P(T)$.
\end{itemize}
\end{example}

\begin{definition}
A map between non-$\Sigma$ operads is a natural transformation which respects
the decomposition of axiom~(\ref{1}) of \refD{op}.
An operad with multiplication is a non-$\Sigma$ operad $\Op$ equipped
with a map from the associative operad $\A$.  
\end{definition}

The notion of operad with multiplication is due to Gerstenhaber
and Voronov \cite{GeVo94}.
The canonical example is the
endomorphism operad of an associative algebra $End(A)$.
Algebras over an operad with multiplication
are in particular associative algebras.
An operad with multiplication in the 
category of spaces is an operad in the category of pointed spaces.

While we have 
taken a categorical approach to defining  operads, 
we will take a more coordinatized approach to their
associated cosimplicial objects.  Recall  
the $\circ_i$ operations $\circ_i : \Op(n) \odot \Op(m) \to \Op(n+m-1)$
which provide a  basic set of morphisms for an operad,
as illustrated in \refF{basicmor}.
From section~3 of \cite{McSm01} we have the following.

\begin{definition} \label{D:hoch}
\begin{itemize}
\item Given an operad with multiplication $\Op$, let 
$\mu$ denote the  morphism $\A(2) = \one_\mC \to \Op(2)$.  
\item Define $d^i : \Op(n) \to \Op(n+1)$ by
$$
d^i = 
\begin{cases}
\one_\mC  \odot \Op(n) \overset{\mu \odot id} \to \Op(2) \odot \Op(n) \overset{\circ_1}{\to}
\Op(n+1) \;\; {\text{if}} \;\; i=0 \\
\Op(n) \odot \one_\mC  \overset{id \odot \mu} \to \Op(n) \odot \Op(2) \overset{\circ_i}{\to}
\Op(n+1) \;\; {\text{if}} \;\;  0<i<n+1\\
\one_\mC \odot \Op(n) \overset{\mu \odot id} \to \Op(2) \odot \Op(n) \overset{\circ_{2}}{\to}
\Op(n+1) \;\; {\text{if}} \;\; i=n+1.
\end{cases}
$$
\item Define $s^i$ as $\Op(c_i)$ where $c_i : \gamma_n \to \gamma_{n-1}$ contracts
the $i$th leaf of $\gamma_n$.
\item Let $\Op^\bullet$ be the cosimplicial object in $\mC$ whose $n$th entry
is $\Op(n)$ and whose coface and codegeneracy maps are given by $d^i$ and
$s^i$ above.  If $\mC$ is the category of vector spaces over a given field, 
let $HH^*(\Op)$ be the homology of the cochain complex defined
by the cosimplicial vector space $\Op^\bullet$.
If $\mC$ is the category of spaces, we call $\tot(\Op^\bullet)$ the 
totalization of $\Op^\bullet$.  
\end{itemize}
\end{definition}

It is straightforward to show that the maps $d^i$ and $s^i$ satisfy
cosimplicial identities.

\begin{remark}
In the category of vector spaces, Tourtchine introduced the terminology $HH^*(\Op)$
because if $A$ denotes an associative algebra and ${\rm End}(A)$ is its
endomorphism operad then $HH^*({\rm End}(A)) = HH^*(A)$, the usual
Hochschild cohomology of $A$. We are not aware, however, of any sense in which
Hochschild  cohomology of operads is a cohomology theory for operads.
Instead, Kontsevich conjectures that there is a suitable enriched homotopy structure
on the category of operads of spaces such that $\tot(\Op^\bullet)$ is the derived
space of maps from the associative operad to $\Op$.
\end{remark}

\section{The choose-two operad}

At the combinatorial heart of our work is the choose-two operad.
Its operad structure intertwines the sets $\binom{\n}{2}$, of distinct pairs of
elements $(i, j) \in \n= \{1,\ldots, n\}$  with $i < j$ for definiteness, with
rooted planar trees.  
Recall that $\Fo$, the opposite category to the category of pointed
sets, is symmetric monoidal with product  given by pointed union,
denoted $\vee$,
and unit given by the one-point set.
Let $S_+$ denote the union of a set $S$ with a disjoint base point.

\begin{definition} \label{D:mB}
\begin{itemize}

\item The join of two leaves in an rp-tree is the first vertex
(that is, the farthest from the root) at which their
root paths coincide.

\item Label both the leaves of an rp-tree and the edges which
emanate from a given vertex $v$ with elements of
$\n$ and $\{1, \ldots, |v|\}$ respectively
according to the order given by  planar embedding.  
To an rp-tree $T$ with $n$ leaves and
two distinct integers $i,j \in \n$ let $v$ be the join
of the leaves labelled $i$ and $j$ and define
$J_v(i), J_v(j)$ to be the labels of the edges of $v$ over
which leaves $i$ and $j$ lie, as illustrated below.

\begin{center}
\begin{mydiagram}\label{F:binomop}

\begin{center}$$\includegraphics[width=10cm]{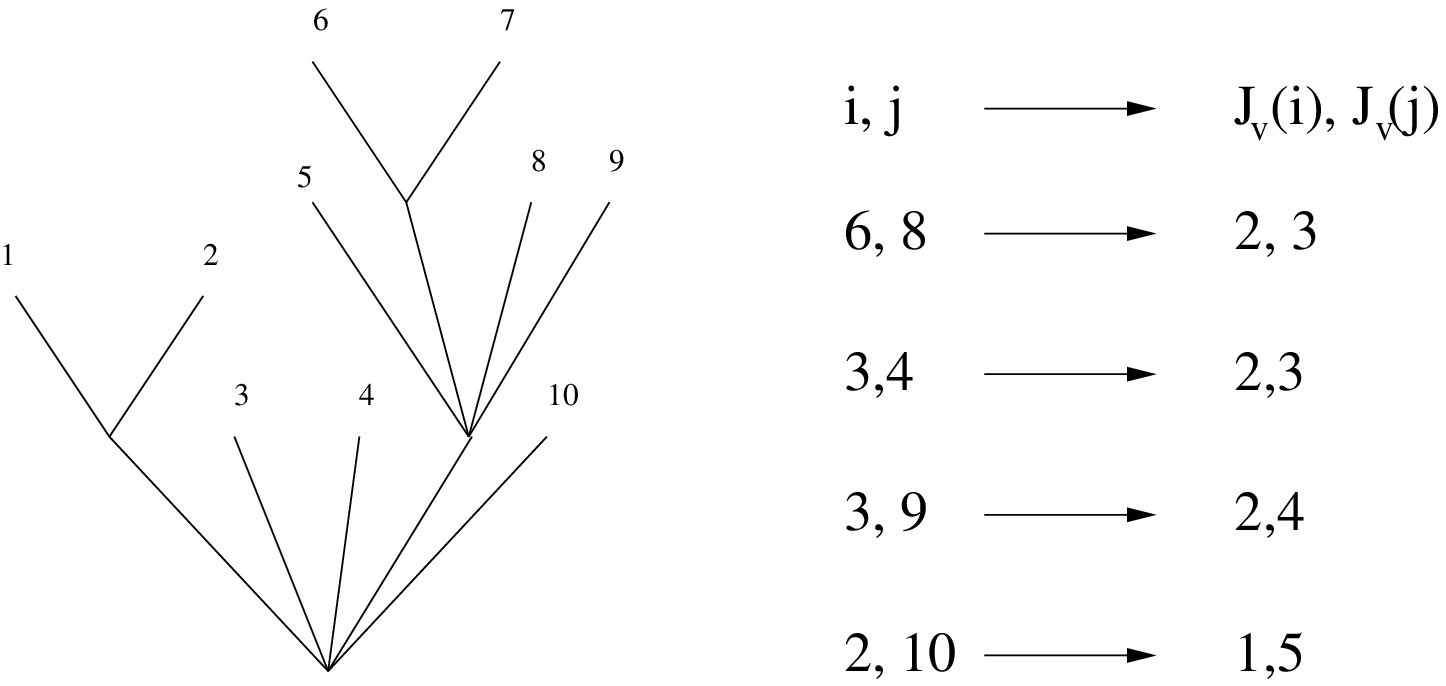}$$
\end{center}
\end{mydiagram}

\end{center}

\item  Let  $\mB$, the choose-two operad, be the non-$\Sigma$ operad in the category
$\Fo$ defined as follows:
\begin{itemize}
\item $\mB(T) = \bigvee_{w} \binom{|w|}{2}_+,$ where $w$ ranges over vertices of $T$.
\item $\mB(T \to \gamma_n)$, where $\gamma_n$ is
a corolla, is the function 
$$(i,j) \in \binom{\n}{2} \mapsto 
(J_v(i), J_v(j)) \in \binom{|v|}{2} \subset \bigvee_{w \in T} \binom{|w|}{2}_+,$$
where $v$ is the join of leaves $i$ and $j$.
\end{itemize}
\end{itemize}
\end{definition}

With our choice of definitions, it is straightforward to verify that $\mB$ is an operad.

\begin{theorem}\label{T:binoms2}
As a simplicial set, 
the cosimplicial object in $\Fo$ associated to the choose-two operad 
$\mB^\bullet$ is
isomorphic to $S^2_\bullet$, the simplicial model for $S^2$.
\end{theorem}

\begin{proof}
Recall that $S^2_\bullet = \Delta^2_\bullet/ \partial
\Delta^2_\bullet$, where $\Delta^2_\bullet$ is the standard 
simplicial model for $\Delta^2$. 
The set $n$-simplices of $\Delta^2$ is the 
set of $(x_0 \leq x_1 \leq \cdots \leq x_n) \in \{0, 1, 2\}^{n+1}$,
so the cardinality of $\Delta^2_n$ is the $(n+1)$st
triangular number.  The $i$th face and degeneracy maps are 
defined by deleting and repeating $x_i$,
respectively.  To obtain $S^2_\bullet$ we identify all $n$-tuples
in which one of $\{0, 1, 2\}$ does not appear to a single simplex
in each degree, which is degenerate in positive degrees. 

The $n$th entry of
$S^2_\bullet$ is isomorphic to 
$\binom{\n}{2}_+$, the set of unordered pairs of points in 
$\n$, along with a disjoint point $+$ which is the image
of $\partial \Delta^2_\bullet$ under the quotient map.
The isomorphism records the indices $j$ and $k$ for which
$x_{j-1} < x_j$ and $x_{k-1} < x_k$, when there are two such 
indices.  When there are not two such indices, such a sequence
is identified with the degenerate point $+$.
Under this isomorphism  $d_i$ sends $+$ to $+$ and for $i \neq 0, n$ sends
\begin{equation}\label{E:1}
(j,k) \mapsto 
\begin{cases}
(\delta_i(j), \delta_i(k)) \;\;  {\rm{if}} \;\;  \delta_i(j) \neq \delta_i(k)\\
+ \;\; {\rm{otherwise}}
\end{cases}
\;\;\;  {\rm{where}} \;\;\; 
\delta_i(j) = 
\begin{cases}
j \;\; {\rm{if}} \;\; j \leq i\\
j-1  \;\; {\rm{if}} \;\; j > i.
\end{cases}
\end{equation}
For $i = 0$ and $i=n$
the basic formula is the same, but the $(j,k)$ which get sent to
$+$ are those with $j=1$ or $k=n$, respectively.  Similarly, $s_i$
sends $+$ to $+$ and sends
$(j,k) \mapsto  (\sigma_i(j), \sigma_i(k))$
where $\sigma_i(j) = j$ if $j \leq i$ or $j+1$ otherwise.

By definition $\mB^n = \binom{\n}{2}_+$, and it is 
straightforward to check that the structure maps of $S^2_\bullet$ and
the associated cosimplicial object of $\mB^\bullet$ coincide.
To give an example, we unravel the definition of $d^i$ with $0<i<n$
for $\mB^n$.
These coface maps are given by composites
$\left( \mB(\gamma_n) \vee +\right) \overset{id \vee \mu}{\longrightarrow} 
\mB(\includegraphics[width=0.6cm]{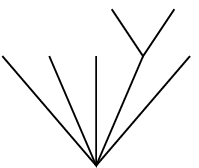}^i_{n}) 
\overset{\circ_i}{\longrightarrow} \mB(\gamma_{n+1})$,
where $\includegraphics[width=0.6cm]{Figures/sm1.eps}^i_{n}$ 
denotes the tree with
$n$ root edges and one trivalent internal vertex, 
which is terminal for the $i$th root edge.
In $\Fo$ the morphism $id \vee \mu$ corresponds to the collapse
map in $\Set$ which sends 
$\binom{\mathbf{2}}{2}_+ \subset \binom{\n}{2}_+ \vee
 \binom{\mathbf{2}}{2}_+ $ to the base point and is the identity on
 $ \binom{\n}{2}_+$.  The morphism $\circ_i$ sends $(i,i+1)$ to 
 $(1,2) \in \binom{\mathbf{2}}{2}$
 and sends all other $(j,k)$ to $(\delta_i(j), \delta_i(k)) \in \binom{\n}{2}$.
 The composite of these two maps coincides with the definition
 of $d_i$ for $S^2_\bullet$, as in \refE{1}.
\end{proof}

Let $X^{\mB^\bullet}$ be the operad which is 
the composite of the operad $\mB^\bullet : \Upsilon \to \Fo$ with 
the symmetric monoidal  functor $X^{-} : \Fo \to \Top$.
 \refT{binoms2} implies the following.

\begin{corollary}
For any $X$ in a category $\mC$ whose categorical product is 
symmetric monoidal,  $X^{S^2_\bullet}$
canonically defines an operad  through its isomorphism
with $X^{\mB^\bullet}$.
\end{corollary}

We have yet to find any familiar interpretation for algebras over this operad
in the categories of spaces and vector spaces.  For spaces the operad
structure on $X^{S^2_\bullet}$ does have consequences, as we explain in
\refX{omega2}.  

\section{The Kontsevich operad}

In this section we define an operad structure on the completion
of configurations in Euclidean space up to translation and scaling defined by
Kontsevich \cite{Kont99}.  The fact that one could define operads using the 
canonical completion of
configuration spaces was noticed by Getzler and Jones \cite{GeJo94}
soon after this completion was introduced by Fulton-MacPherson \cite{FuMa93} and 
Axelrod-Singer \cite{AxSi94}.  This operad structure was fully developed by 
Markl \cite{Mark99}.   The variant with which we work was first proposed by
Kontsevich \cite{Kont99}, but Gaiffi \cite{Gaif03}
pointed out the deviation with the canonical completion.   Indeed, while Kontsevich
called the following the Fulton-MacPherson operad, we call it the
Kontsevich operad to highlight the difference between the two.
Though this construction lacks some of the properties of the
canoncial completion, in particular smoothness, is has diagonal and
projection maps which satisfy simplicial identities exactly rather than
up to homotopy.  These properties led to this construction's independent discovery,
its use, and its naming as the simplicial variant in \cite{Sinh02}.

We start by setting notation for products of spaces and maps, which we will use
extensively.

\begin{notation}
If $S$ is a finite set, let $X^S$ be the space of all functions from $S$
to $X$ with the product topology.  
For coordinates we use
$(x_s)_{s\in S}$ or just $(x_s)$ when $S$ is understood. Similarly,  a
product of maps $\prod_{s\in S} f_s$ may be written$(f_s)_{s \in S}$ or
just $(f_s)$.  Recall that $\n = \{ 1, \ldots, n \}$.
\end{notation}

\begin{definition}\label{D:cnrm}
\begin{itemize}
\item Let $C_n(\R^m)$ denote the space of $(x_i) \in (\R^m)^{\n}$ such that 
if $i \neq j$ then $x_i \neq x_j$. 
Let $\wt{C}_n(\R^m)$ be the quotient of  $C_n(\R^m)$ by the equivalence 
relation generated by translating all of the $x_i$ by some $v$
or multiplying them all by the same positive scalar.
\item For any $v \in \R^m - 0$, let $u(v) = \frac{v}{||v||}$, the unit
vector in the direction of $v$.
\item Let $\wt{C}_n \la \R^m \ra$
be the closure of the image of $\wt{C}_n(\R^m)$
under the map $(\pi_{ij})$ to $(S^{m-1})^{\binom{\n}{2}}$, where 
$\pi_{ij}$ sends the equivalence class of $(x_i)$ to $u(x_i - x_j)$.
\end{itemize}
\end{definition}

Note that $(\pi_{ij})$ is not injective -- it fails to be so on configurations
in which all the $x_i$ lie on some line
-- so $\wt{C}_n (\R^m)$ is  not a subspace of $\wt{C}_n \la \R^m \ra$.
But we do have the following theorem, a consequence of 
Corollary~4.5, Lemma~4.12 and Corollary~5.10 of \cite{Sinh03}.

\begin{theorem}\label{T:inclheq}	
The canonical map $\wt{C}_n (\R^m) \to \wt{C}_n \la \R^m \ra$ 
is a homotopy equivalence.
\end{theorem}

What makes $\wt{C}_n \la \R^m \ra$ manageable is that we can characterize
it as a subspace of $(S^{m-1})^{\binom{\n}{2}}$.   We extend
coordinates for $(u_{ij}) \in (S^{m-1})^{\binom{\n}{2}}$ by letting $u_{ji}$ be $-u_{ij}$
when $j>i$. 

\begin{definition}
\begin{itemize}

\item A chain, or $k$-chain, in $S$ is a collection 
$\{ i_1 i_2,$ $i_2 i_3$, \ldots, $i_{k-1} i_k \}$, with all
$i_j \in S$ and $i_j \neq i_{j+1}$.  
Such indices label the edges of a path in the complete
graph on $S$.  A chain is a loop, or $k$-loop,  if $i_k = i_1$.
A chain is straight if  it does not contain any loops.
The reversal of a chain is the chain $i_k i_{k-1}, \ldots, i_2 i_1$.

\item  A point $(u_{ij}) \in (S^{m-1})^{\binom{\n}{2}}$ is three-dependent 
if for any $3$-loop
 $L$  in $\n$ there exist $a_{ij} \geq 0$,
with at least one non-zero, such that $\sum_{ij \in L} a_{ij} u_{ij} = 0$.

\item  If $S$ has cardinality four and is ordered, 
we may associate to a straight $3$-chain $C$
a permutation of $S$ denoted $\sigma(C)$ which relates 
the order in which 
indices appear in $C$ to the ordering on $S$.  For example, 
if $S = \{1, 2, 3, 4\}$, then $\sigma(23, 31, 14) = (2314)$.
A complementary $3$-chain $C^*$ is a chain, unique up
to reversal, which comprises the three pairs of indices not in $C$.

\item A point $(u_{ij}) \in (S^{m-1})^{\binom{\n}{2}}$ is four-consistent if for any $S \subset \n$ of 
cardinality four and any $v, w \in S^{m-1}$
we have that 
\begin{equation} \label{4consist}
\sum_{C \in \mathcal{C}^3(S)} (-1)^{|\sigma(C)|}
\left( \prod_{ij \in C} u_{ij} \cdot v \right)\left(\prod_{ij \in C^*} u_{ij} \cdot w \right) = 0,
\end{equation}
where $\mathcal{C}^3{S}$ is the set of straight $3$-chains in $S$
modulo reversal and ${|\sigma(C)|}$ is the parity of $\sigma(C)$.
\end{itemize} 
\end{definition}

Points in the image of $C_n(\R^m)$ under $(\pi_{ij})$ are 
three-dependent and four-consistent, and also satisfy $u_{ij} = -u_{ji}$, 
a condition we refer to as anti-symmetry.
These properties also hold for $C_n \la \R^m \ra$, the closure,
by continuity.  Adding the converse, 
we have the following, which is Theorem~5.14 of \cite{Sinh03}.

\begin{theorem}\label{T:cnlasub}
$\wt{C}_n\la\R^m\ra$ is the subspace of all three-dependent, four-consistent points
in $(S^{m-1})^{\binom{\n}{2}}$.
\end{theorem}

We will define operad maps on the completions
$\wt{C}_n \la \R^m \ra$ through 
coordinates of $(S^{m-1})^{\binom{\n}{2}}$.  Embed
$(S^{m-1})^{\binom{\n}{2}}$, and thus $\wt{C}_n \la \R^m \ra$, in 
$(S^{m-1})^{\binom{\n}{2}_+}$ as the subspace of $(u_{ij}) \times u_+$ with
$u_+$ equal to the basepoint of $S^{m-1}$, which we choose to be
the south pole $*_S = (0, \ldots, 0, -1)$.

\begin{theorem} \label{T:kop}
The operad structure on $(S^{m-1})^{\mB^\bullet}$ restricts to the 
subspaces $\wt{C}_n \la \R^m \ra$. 
\end{theorem}

We call the resulting operad, whose $n$th entry is $\wt{C}_n \la \R^m \ra$,
the Kontsevich operad $\K_m$.

\begin{proof}
Given a tree $T$, let $(u^v_{k\ell})$ be a point in $(S^{m-1})^{\mB(T)}$,
where $v$ ranges over vertices of $T$ and $k,\ell \in \binom{{\mathbf{|v|}}}{2}$.
By \refD{mB}, the operad structure on $(S^{m-1})^{\mB^\bullet}$
sends the morphism $T \to \gamma_n$ to the map given in coordinates by 
$(u^v_{k \ell}) \mapsto (w_{i j})_{i,j \in \binom{\n}{2}}$, where $w_{i j} = u^v_{J_v(i), J_v(j)}$ and 
$v$ is the join vertex of the leaves $i$ and $j$.  

We verify that if the $(u^v_{k \ell})$ satisfy three-dependence and
four-consistency for each $v$, then so does $(w_{ij})$.  For three-dependence,
given some $w_{ij}$, $w_{jk}$ and $w_{ki}$, there are two cases
to consider.  In the first case the join in $T$ of 
leaves $i$ and $j$ lies over that of $i$ and $k$,
so that $w_{jk} = -w_{ki}$ or $0 w_{ij} + 1 w_{jk} + 1 w_{ki} = 0$.
In the second case the joins of $i$ and $j$ and $k$ are all equal to the
same $v$, so that the dependence of $w_{ij}$, $w_{jk}$ and
$w_{ki}$  follows from that of $(u^v_{J_v(i) J_v(j)})$, $(u^v_{J_v(j) J_v(k)})$ 
and $(u^v_{J_v(k) J_v(i)})$.   Four-consistency works similarly.  Given indices
$i, j, k$ and $\ell$ the pairwise joins could all equal some $v$, in which case four
consistency of these $\{ w_{ij} \}$ follows from that of $\{ u^v_{k \ell} \}$.  Or,
if for example the join of $i$ and $j$ lies over those of $i$, $k$ and $\ell$,
then $w_{ik} = w_{jk}$ and $w_{i \ell} = w_{j \ell}$, so four-consistency will
follow by the canceling of terms which agree but for opposite signs.
\end{proof}

In \cite{Sinh03}, we stratify $\wt{C}_n \la \R^m \ra$, and in particular
the points added in closure. 
We will not need this stratification explicitly for our applications, but the related
geometry is helpful in  
understanding the operad structure of $\K_m$.  The stratification is indexed by
rp-trees with no redundant edges, with the $T$th stratum being the image of a map from 
$\prod_{v \in V(T)} \wt{C_{|v|}}(\R^m)$ to $\wt{C}_n \la \R^m \ra$ sending
$(x^v_i) \mapsto (u_{ij})$  with $u_{ij} = \pi_{J_v(i) J_v(j)}((x_i)_v)$. 
Studying this stratification led us to the definition of the choose-two operad.
See Section~3 and Theorem~5.14
of \cite{Sinh03} for a full development of this geometry, which is illustrated
in \refF{cnop}.

\begin{center}
\begin{mydiagram}\label{F:cnop}
The effect of an operad structure map associated 
to the morphism $\includegraphics[width=2cm]{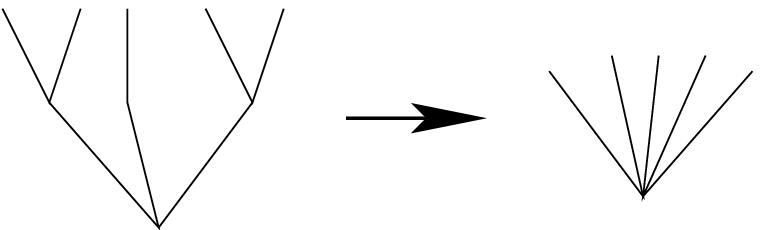}$.
\begin{center}
$$\includegraphics[width=12cm]{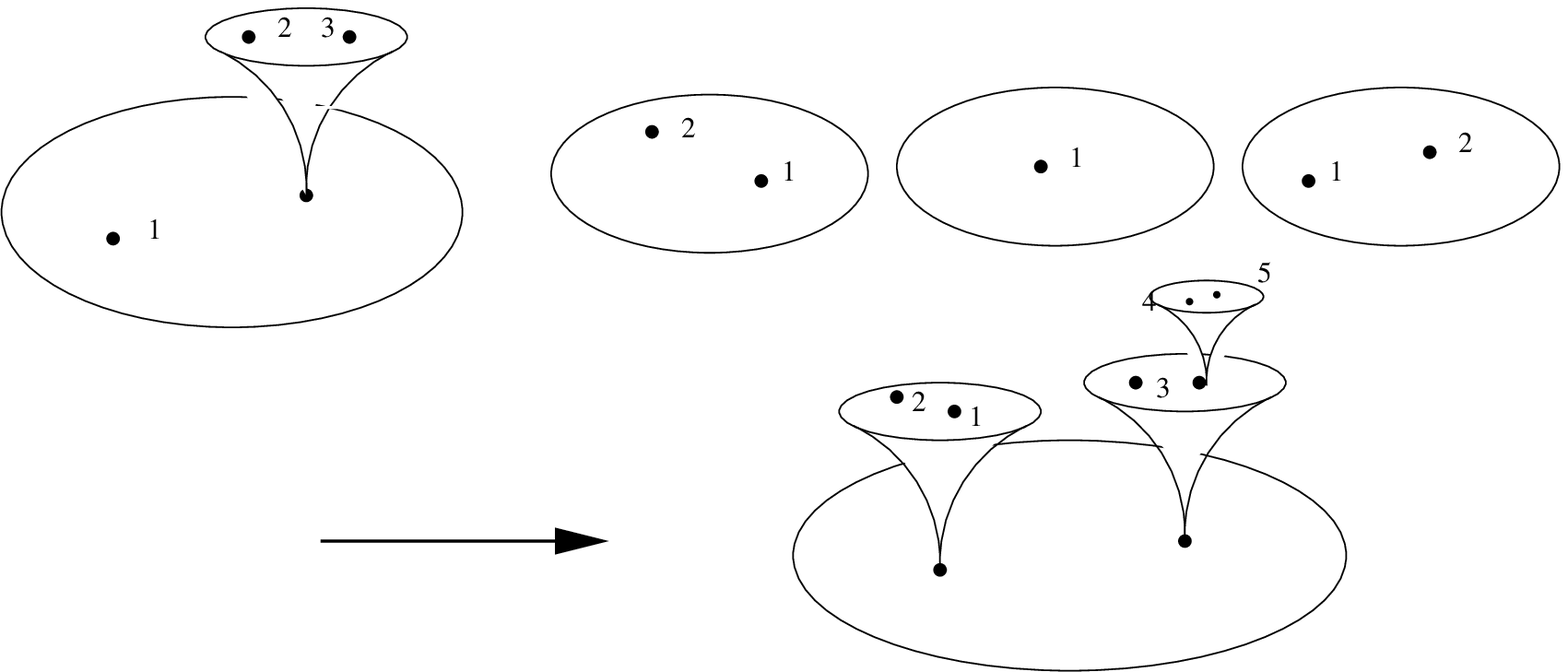}$$
\end{center}
\end{mydiagram}

\end{center}

The standard completions $\wt{C_n}[\R^m]$ also constitute entries of
an operad, which 
has been more intensively studied \cite{GeJo94, Kont99, Mark99}.  The 
reason we use $\K^m$ is the following.   

\begin{proposition}\label{P:assoctoK}
The associative operad maps to the Kontsevich operad, for definiteness
by choosing the basepoint $(x_{ij}) \in \wt{C_n} \la \R^m \ra$ with all $x_{ij} = *_S$,
for all $n$.
\end{proposition}

Finally, we give a comparison between the little disks operad, which we
need to formalize, and the Kontsevich operad.  

\begin{definition}\label{D:littledisk}
\begin{itemize}
\item Recall that the space of $n$ little disks in $D^m$,
the unit disk, denoted $D^m(n)$ is the subspace of $C_n(D^m) \times (0,1]^n$ of 
$(x_i) \times (r_i)$ such that the balls $B(x_i, r_i)$ are contained in 
$D^m$ and have disjoint interiors.  

\item Let $T$ be a tree whose vertices consist of the root
vertex $v_0$ and a terminal vertex $v_e$ for each root edge $e$.  
Thus, $T \to \gamma_n$, where $n$ is the number of leaves of $T$,
gives rise to one of May's structure maps as in \refF{basicmor}.  
Given a label $i \in \n$ let $v(i)$ be the initial vertex for the $i$th
leaf, let $o(i)$ be the label of leaf $i$ within the ordering on edges
of $v(i)$ and let $e(i)$ be the label of the root edge for
which $v(i)$ is terminal.

\item Define $\D^m(T \to \gamma_n)$ as follows 
$$(x^v_i, r^v_i)^{v \in V(T)}_{1 \leq i \leq \#v} \mapsto (y_j, \rho_j)_{j \in \n} \;\;\; {\rm where}
\;\;\; y_j = x^{v_0}_{e(j)} + r^{v_0}_{e(j)}  x^{v(j)}_{o(j)} \;\;\; 
{\rm and} \;\;\; \rho_j = r^{v_0}_{e(j)} r^{v(j)}_{o(j)}.$$
\end{itemize}
\end{definition}

Boardman and Vogt  \cite{BoVo73} and May \cite{May72} showed that 
connected algebras over $\D^m$ are $m$-fold loop spaces.

\begin{theorem}\label{T:DKcomp}
Let $T$ be a tree with a vertex over each root edge as in \refD{littledisk}
above.  The following diagram
commutes up to homotopy,
$$
\begin{CD}
\D^m(T) @>{\D^m(T \to \gamma_n)}>> \D^m(n) \\
@V{p_T}VV @V{p_n}VV\\
\K^m(T) @>{\K^m(T \to \gamma_n)}>> \K^m(n), 
\end{CD}
$$
where the vertical maps $p_T$ are the
products of projections $p_n : D^m(n) \to C_n (\R^m)$ composed with the canonical
maps $C_n(\R^m) \to \wt{C}_n(\R^m) \to \wt{C}_n \la \R^m \ra$. 
Moreover, the vertical maps are homotopy equivalences.
\end{theorem}

\begin{proof}
We define the homotopy explicitly.  
Define $H : \D^m(T) \times (0,1] \to \wt{C}_n(\R^m)$ by sending
$(x^v_i, r^v_i)$ as in \refD{littledisk} and $t \in (0,1]$ to the 
equivalence class of 
$(y_j(t))$ with $y_j(t) = x^{v_0}_{e(j)} + t \cdot r^{v_0}_{e(j)}  x^{v(j)}_{o(j)}$.   
We claim that $H$ extends uniquely to  
$\overline{H} : \D^m(T) \times [0,1] \to \wt{C}_n \la \R^m \ra$, and
 that $\overline{H}$ coincides with 
$\K^m(T \to \gamma_n) \circ p_T$ when $t = 1$.  
Consider $u_{ij} = \pi_{ij}\left( (y_k(t)) \right)$.  If the join of leaves
$i$ and $j$ is one of the non-root vertices,
so  $v(i) = v(j)$, then $u_{ij}$ will be equal to the unit 
vector in the direction of $x^{v(j)}_{o(j)} - x^{v(j)}_{o(i)}$,
independent of $t$.  If the join of leaves $i$ and $j$ is the 
root vertex, then as $t$ approaches $0$, $u_{ij}$
approaches the unit vector in the direction of 
$x^{v_0}_{e(j)} - x^{v_0}_{e(i)}$.  These limiting
values coincide with the definition of $\K^m(T \to \gamma_n) (x^v_i)$.

That the projection $D^m(n) \to C_n(D^m)$
is a homotopy equivalence is standard, known since the 
definition of little cubes in \cite{Boar71}, so by \refT{inclheq} the 
maps $p_n$ are homotopy equivalences.
\end{proof}

In fact, \cite{Kont99} claims that these $\D^m$ and $\K^m$ 
are homotopy equivalent operads, which we assume to 
mean that there is a chain of equivalences of maps of operads, 
that is maps which commute with structure maps exactly.  
We will not need this stronger claim.  
Recall that the homology of an operad of spaces with field 
coefficients is an operad of vector spaces by the K\"unneth theorem.
The homology of the little disks operad has a well-known description.

\begin{definition}\label{D:poiss}
The $k$th entry of the degree $n$
Poisson operad ${\rm Poiss}_n(k)$ is the submodule of the
symmetric algebra on the free graded Lie algebra over $k$ variables
$x_1, \ldots, x_k$  spanned by 
monomials in which all variables appear exactly once.  Monomials
are graded by putting all $x_i$ in degree zero and giving the
bracket degree $n$.   So for example $[x_1, x_3] \left[ [x_4, x_2] x_5 \right]$
and $x_1x_2 \cdots x_5$ are elements of ${\rm Poiss}_3(5)$ of degree
nine and zero respectively.

The map $\circ_i : {\rm Poiss}_n(j) \otimes {\rm Poiss}_n(k) \to {\rm Poiss}_n(j+k-1)$ 
sends $m_1 \otimes m_2$ to the monomial defined as follows.
\begin{itemize}
\item For each $j$, substitute $x_{j+i-1}$ for $x_j$ in $m_2$ to obtain
$\overline{m_2}$.
\item In $m_1$, substitute $x_{j+i-1}$ for $x_j$ if $j > i$ and 
$\overline{m_2}$ for $x_i$  to obtain $m$.
\item Reduce $m$ according to the graded Leibniz rule
$$ [a, bc] = [a,b] c + (-1)^{(|a| + n + 1)|b|} b [a,c],$$
to obtain an element of ${\rm Poiss}_n(i+j-1)$.
\end{itemize} 
\end{definition}

The following corollary is essentially a summary of 
Fred Cohen's famous  calculation of the homology of $\D^m$
\cite{FCoh76}.   We also plan to give an exposition of this
result in \cite{Sinh04.2}.

\begin{theorem}\label{T:poisson}
The homology of $\K^m$ is ${\rm Poiss}_{m-1}$.
\end{theorem}


\section{Models for spaces of knots and immersions arising from the calculus of embeddings}

\subsection{A brief overview of the calculus of embeddings}

Our main theorems connect the theory of operads to 
Goodwillie calculus.  We first informally introduce some terminology
from the calculus of embeddings (see Weiss's \cite{Weis95} for an 
excellent introduction and \cite{Weis99} for a full treatment), and 
then precisely state the theorems we use.  The main spaces
with which we are concerned are related to embeddings and immersions.

\begin{definition}\label{D:embimm}
Let $\Emb(M,N)$ denote the space of embeddings of $M$ in $N$, topologized
as a subspace of the space of $C^1$ maps.  
Similarly, let $\Imm(M,N)$ be the space of immersions of $M$ in $N$.  
If $M$ and $N$ have boundary, we usually specify some boundary conditions.  In
particular, if $M = \I$ and $N = \I^m$, we let $*_+  = (0, \ldots, 0, 1) \in \I^m $, 
$*_- = (0, \ldots, 0, -1)$ and demand that the endpoints of $\I$ map to $*_+$ and $*_-$
with tangent vectors $*_S$.
\end{definition}

By results of Palais \cite{Pala66}, these spaces are dominated by
simplicial complexes and thus homotopy equivalent to
CW-complexes \cite{Misl95}.

In the calculus of embeddings, we view spaces of embeddings, 
immersions and other moduli in differential topology as functors from the
poset  of open subsets of $M$ to topological spaces, a philosophy
originally due to Gromov.  Ultimately interested in the value of the functor at the open
set which is all of $M$, we try to use homotopy limits to
interpolate that value from 
values of the functor at  simple open sets, namely those which
are diffeomorphic to a union of finitely many disjoint open balls.  Functors for which 
interpolation using a finite number of balls works perfectly are called 
polynomial, and those for which interpolation works in the limit as the number
of balls tends to infinity are 
called analytic.  Weiss shows in \cite{Weis99} that polynomial 
functors are those which satisfy higher-order Mayer-Vietoris 
conditions, and Goodwillie-Weiss show in \cite{GoWe99} 
that analyticity follows from satisfying those 
conditions through an increasing range of connectivity.
More formally we have the following.

\begin{definition}\label{D:Tk}
\begin{itemize}
\item For any manifold $W$ of dimension $m$ 
let $\U(M)$ be the category of open subsets of $M$ under
inclusion, and let $\U_k(M)$ be the full sub-category of $\U(M)$
of open sets diffeomorphic to $\sqcup_i \R^m$, where
$i \leq k$. 
\item For any contravariant functor $F$ from $\U(M)$ to 
spaces let $T_k F$ be the functor which sends $W$
to $\holim{U \in \U_k(W)}{F(U)}$.
\item Let $b_k(F) : T_k F \to T_{k-1} F$ be the canonical natural
transformation defined by restricting $\U_k(M)$ to $\U_{k-1}(M)$,
and let $T_\infty F$ be the homotopy inverse limit of the $T_k F$
over these restrictions.
\item Let $\eta_k(F) : F \to T_k F$ be the canonical natural
transformation arising from the maps $F(M) \to F(U)$
for $U \in \U_k(M)$.  If by context $F$ is understood, we will
use the simpler notation $\eta_k$.
\item The natural transformations $\eta_k$ commute with the $b_k$,
so let $\eta_\infty : F \to T_\infty F$ be the limiting natural
transformation.
\end{itemize} 
\end{definition}

The sequence $T_0 F \overset{b_1}{\leftarrow} 
T_1 F \overset{b_2}{\leftarrow} T_2 F \leftarrow \cdots$ is called the
Taylor tower for $F$.  Analyticity means that the homotopy
inverse limit of this tower is weakly equivalent to $F$.
The motivating example for this circle of ideas is that
of immersions.

\begin{theorem}\label{T:imm}
If $\dim(M) < \dim(N)$ then for $k \geq 1$, 
$\eta_k : \Imm(U,N) \to T_k \Imm(U, N)$
is a weak equivalence for any $U \subseteq M$.
\end{theorem}

This theorem follows from Example~2.3 of \cite{Weis99}, which says 
that $\Imm(-, N)$ is a linear functor, and the commentary 
after Theorem~5.1 of \cite{Weis99}.   
See also \refT{modimm} below and Proposition~5.12 of \cite{Voli04.1}.
The embedding functor is not polynomial but by theorems
of Goodwillie, Klein and Weiss it is analytic.  The following Theorem is
essentially Corollary~2.5 of \cite{GoWe99}. 

\begin{theorem}\label{T:GK}
If $\dim(M) < \dim(N) -2$ then $\eta_\infty (\Emb)$ is a weak equivalence.
\end{theorem}

For $\dim(M) < \dim(N) -2$ as stated, this theorem requires deep
disjunction results of Goodwillie, and surgery results of Goodwillie-Klein
\cite{Good03, GoKl03}.  If $\dim(M)$ is less than roughly $\frac{\dim(N)}{2}$,
in particular when $M$ is one-dimensional and $N$ has dimension
five or greater, there
are much easier methods, using only the generalized Blakers-Massey theorem
and dimension counting, for proving the needed
higher-order Mayer-Vietoris conditions.

\subsection{Knot space models through homotopy limits of configuration spaces }\label{knotmod}

\refD{Tk} of $T_k \Emb(M,N)$ seems difficult to manage, being a homotopy
limit over a  large indexing category.  But the building blocks, 
namely spaces of embeddings of balls, are essentially configuration
spaces.  Goodwillie, Klein and Weiss
have used this observation to give
more concrete models for the spaces in this Taylor tower
(or for the homotopy fibers of $T_k \to T_{k-1}$, which are called layers), 
either as spaces of sections, as in section nine of \cite{Weis99}, 
or as mapping spaces with strongly defined equivariance properties,
as in \cite{GKW01}.  In the case of knots we have developed three
closely-related models for these polynomial approximations
\cite{Sinh02} and used them for both computational and
geometric applications \cite{ScSi02, BCSS03}.  These models all 
utilize completions of configuration spaces constructed similarly to
the Kontsevich operad.

\begin{definition}
\begin{itemize}
\item Let $A_n\la \I^m \ra$ be the product $(\I^m)^{\n} \times
(S^{m-1})^{\binom{\n}{2}}$, with coordinates $(x_i) \times (u_{ij})$.
\item Let $C_n \la \I^m \ra$ be the closure of the image of $C_n(\I^m)$
under $\iota \times (\pi_{ij})$, where $\iota$ is the inclusion of
$C_n(\I^m)$ in $(\I^m)^{\n}$.  
\item Let $C_n \la \I^m, \partial \ra$ be the closure in $C_{n+2} \la \I^m \ra$
of the subspace of  $C_{n+2}(\I^m)$ with $x_1 = *_+ = (0, \ldots, 0, 1)$ and $x_{n+2} = *_-$.
\end{itemize}
\end{definition}

In \cite{Sinh03} we study $C_n \la \I^m \ra$ by
relating it to the canonical compactification $C_n [\I^m]$, which
is a manifold with corners.  We characterize
$C_n \la \I^m \ra$ as a subspace of its defining ambient space,
as stated for $\wt{C}_n\la \R^k \ra$ in \refT{cnlasub}.  The following
is essentially Theorem~5.14 of \cite{Sinh03}.

\begin{theorem}\label{T:cnlasub2}

$C_n \la I^k \ra$ is the subspace of $(x_i) \times (u_{ij})$ such that 
$(u_{ij}) \in \wt{C}_n \la \R^k \ra$ and if $x_i \neq x_j$ then 
$u_{ij}$ is $u(x_j - x_i)$.
\end{theorem}

In our models, we need diagonal 
maps between configuration spaces.  The idea is to add a point 
``infinitesimally far'' from one point in a configuration,
but to do so entails choosing a unit tangent vector at that point.

\begin{definition}
Let $C_n' \la \I^m \ra = C_n \la \I^m \ra \times (S^{m-1})^{\n}$.
Let $A'_n\la \I^m \ra = (\I^m \times S^{m-1})^{\n} \times
(S^{m-1})^{\binom{\n}{2}}$, which is canonically diffeomorphic to  
$(\I^m)^{\n} \times (S^{m-1})^{\frac{\n(\n + 1)}{2}}$.  We use coordinates for 
this latter presentation of the form
$(x_i) \times (u_{ij})$ with $i \leq j$, and if $i > j$ we
set $u_{ij} = -u_{ji}$.
\end{definition}

As  in \refT{kop} we define maps between the
$C_n' \la \I^m \ra$ at the level of the ambient spaces
$A'_n\la \I^m \ra$, using \refT{cnlasub2} to check that
they restrict appropriately.  
We are aided by the following combinatorial shorthand.

\begin{definition}
\begin{itemize}
\item Given a map of sets $\sigma : R \to S$ let
$p^X_{\sigma}$, or just
$p_{\sigma}$, denote the map from $X^S$ to $X^R$ which sends 
$(x_i)_{i \in S}$ to $(x_{\sigma(j)})_{j \in R}$.

\item Given $\sigma : \mathbf{m} \to \n$, 
define $A_\sigma : A'_n\la \I^m \ra \to A'_m \la \I^m \ra$ as 
$p^{\I^m}_\sigma \times p^{S^{k-1}}_{\sigma^{(2)}}$, where 
$\sigma^{(2)} = \sigma \times \sigma|_{\binom{\n}{2}}$.
\end{itemize}
\end{definition}

\begin{proposition}[Proposition 6.6 of \cite{Sinh03}]
The restriciton $A_\sigma$ to $C'_n\la \I^m \ra$ maps to $C'_m \la \I^m \ra$.  
If $\sigma$ sends $1 \to 1$
and $n \to m$ then $A_\sigma$ also restricts to a map,
which we call $F_\sigma$, from $C'_{n-2} \la \I, \partial \ra$
to $C'_{m-2} \la \I, \partial \ra$.
\end{proposition}

We may now define diagonal maps on compactified configuration spaces
with tangential data. 

\begin{definition}
Let $\delta^i : C'_n \la \I^m, \partial  \ra \to C'_{n+1} \la \I^m, \partial  \ra$ 
be $F_{\sigma_i}$ where $\sigma_i : \und{n+3} \to \und{n+2}$ sends $j$ to itself 
if $j \leq i$ or $j-1$ if $j>i$.  
\end{definition}

A final key property of this compactification is that it is functorial
for embeddings.  The proof of the following theorem is identical to
that of Corollary~4.8 of \cite{Sinh03}, using Theorem~5.8 of \cite{Sinh03}
and the fact that $C_n \la \I \ra = \Delta^n$.  Recall that for a nonzero
vector $v\in \R^m$, $u(v) = \frac{v}{||v||}$.

\begin{theorem}
For an embedding $f : \I \to \I^k$ there is an evaluation map
$ev_n(f) : \Delta^n \to C_n \la \I^k \ra$ which extends the 
map from the interior of $\Delta^n$ to $C'_n (\I^k)$ sending
$(t_i)$ to $(f(t_i)) \times (u(f'(t_i)))$.
\end{theorem}

One of the main themes of \cite{Sinh02} and of 
\cite{Voli04.1} is connecting this
evaluation map with the calculus of embeddings.
Applying this calculus to embeddings of the unit interval  is 
simpler than to embeddings of higher-dimensional manifolds
because the category $\U_k(\I)$ 
may be replaced by the category of subsets of a finite set
(see \refD{P0}).

\begin{definition} 
Let $\D^m_k$ be the functor from $P_0(k)$ to spaces which sends
$S \subseteq [k]$ to $C'_{\#S - 1} \la \I^m, \partial \ra$ and which sends
the inclusion $S \subset S \cup j$ to the map $\delta^i$ where $i$
is the number of elements of $S$ less than $j$.
Let $D^m_k = \holim{}{ \D^m_k}$.
\end{definition}

In the notation of \cite{Sinh02}, $\D^m_k$ would be $\D_k \la \I^m \ra$.
Because the realization of $P_0(k)$ is $\Delta^k$ and all
of the maps $\delta^i$ are inclusions of subspaces, $D^m_k$
is a subspace of ${\rm Maps}(\Delta^k, C_k'\la \I^m, \partial \ra)$.
If $f \in \Emb(\I, \I^m)$ is a knot, $ev_k(f)$ defines an element
of $D^m_k$, as we may simply check that if $t_j = t_{j+1}$
for some point $(t_i) \in \Delta^k$ then the image of $ev_k(f)\left((t_i)\right)$
is in the image of $\delta^j$.  By abuse, let $ev_k$ denote the adjoint map
from $\Emb(\I, \I^m)$ to $D^m_k$.
Building on the simpler
``cutting method'' definition of  $T_k \Emb(\I, \I^m)$, as described in 
Section~3 of \cite{Sinh02},
Lemma~5.18 and the proof of Theorem~5.3 of \cite{Sinh02}
establish the following.

\begin{theorem}\label{T:modemb}
$D^m_k$  is homotopy equivalent to $T_k \Emb(\I, \I^m)$, and 
$ev_k$ agrees with $\eta_k$ in the homotopy category.  
\end{theorem}

We next take a similar point of view for immersions of an interval,
in order to arrive at our model for $\F_m$. 
This point of view is taken in Section~5.3 of \cite{Voli04.1}, which
also considers the homotopy fiber of the map from knots
to immersions.  Indeed, our \refT{modimm} and \refP{emimsq} below
overlap significantly with Propositions~5.12 and 5.13 of \cite{Voli04.1}.

\begin{definition}
\begin{itemize}
\item Let $d^i : (S^{m-1})^j \to (S^{m-1})^{j+1}$ be the $i$th diagonal inclusion,
which on coordinates repeats the $i$th entry.
By convention, for $i=0$ and $i=j+1$ we insert the basepoint
$*_S$ as the first, respectively last, coordinate.
\item Let $\G^m_k$ be the functor from $P_0(k)$ to spaces which
sends $S \subseteq [k]$ to $(S^{m-1})^{\#S - 1}$ and which sends
the inclusion $S \subset S \cup j$ to the diagonal map $d^i$ where $i$
is the number of elements of $S$ less than $j$.
\item Let $G^m_k = \holim{}{ \G^m_k}$.
\end{itemize}
\end{definition}

As was true for $D^m_k$, $G^m_k$ is a subspace of the space
of maps from $\Delta^k$ to the terminal space of $\G^m_k$, namely
$(S^{m-1})^k$.
The evaluation map for immersions is the unit derivative map.
By abuse, let $ev_k : \Imm(\I, \I^m) \to G^m_k$ send $f$ to the map
which sends $t_1, \ldots, t_k$ to $uf'(t_1), \ldots, uf'(t_k)$.

\begin{theorem}\label{T:modimm}
If $k \geq 1$, $G^m_k$  is homotopy equivalent to $T_k \Imm(\I, \I^m)$,
and thus to $\Imm(\I, \I^m)$.  Moreover, 
$ev_k$ agrees with $\eta_k$ in the homotopy category.  
\end{theorem}

\begin{proof}[Sketch of proof]
There are many ways to establish this theorem.  
By the Hirsch-Smale theorem \cite{Smal59}, $\Imm(\I, \I^m)$
is homotopy equivalent to $\Omega S^{m-1}$,
through the unit derivative map.  But $ev_1 : \Imm(\I, \I^m) \to
\holim{}{ (* \to S^{m-1} \leftarrow *)}$ is also the unit
derivative map, which establishes
the theorem for $k=1$.  For the other $k$,
we may use \refL{coshol}, since $\G^m_k$ is
the pull-back of the standard cosimplicial model for
$\Omega S^{m-1}$ through the functor $c_k$ of \refD{ck}.  The 
$k$th totalization of this cosimplicial model, which is fibrant, 
is homeomorphic to $\Omega S^{m-1}$ if $k \geq 1$,  
from which it follows that $\G^m_k$ is homotopy equivalent
to $\Omega S^{m-1}$.   The map from $\Omega S^{m-1}$
to the $k$th totalization, and thus $G^m_k$, is
through evaluation of the unit derivative.
\end{proof}

Let $\tau: \Emb(\I, \I^m) \to \Imm(\I, \I^m)$ denote the inclusion.
Let $\rho^m_k : \D^m_k \to \G^m_k$ denote the map of diagrams
defined on each entry by projection from 
$C_n' \la \I^m, \partial \ra = C_n \la \I^m, \partial \ra \times (S^{m-1})^n$
onto $(S^{m-1})^n$, and let $p^m_k$ also denote the
induced map on homotopy limits.

\begin{proposition}\label{P:emimsq}
The square
$$
\begin{CD}
\Emb(\I, \I^m) @>{\tau}>> \Imm(\I, \I^m) \\
@V{ev_k}VV             @V{ev_k}VV  \\
D^m_k @>{p^m_k}>>           G^m_k.
\end{CD}
$$
commutes.  Moreover, $p^m_k$ agrees with $T_k(\tau)$ in the homotopy
category.
\end{proposition}

\begin{proof}[Sketch of proof]
The commutativity of the diagram is immediate from the definitions.  That 
$p^m_k$ agrees with $T_k (\tau)$ in the homotopy category ultimately
follows from the fact that for $U$ a disjoint union of $k+2$ open intervals,
two of which contain endpoints of $\I$ and thus are fixed at one end,  we have
$\Emb(U, \I^m) \simeq C_k' \la \I^m, \partial \ra$, $\Imm(U, \I^m) \simeq (S^{m-1})^k$
and the inclusions from embeddings to immersions coincides with projection,
as in the definition of $p^m_k$. 
\end{proof}

\subsection{A closer look at $\F_m$}\label{S:closer}

\begin{proposition}\label{P:null}
The inclusion $\tau: \Emb(\I, \I^m) \to \Imm(\I, \I^m)$ is null-homotopic, so 
$$\F_m  \simeq \Emb(\I, \I^m) \times \Omega \Imm(\I, \I^m) \simeq 
              \Emb(\I, \I^m) \times \Omega^2 S^{m-1}.$$
\end{proposition}

\begin{proof}
Given $f \in \Emb(\I, \I^m)$ consider the map $\rho(f) :  \Delta^2 \to S^{m-1}$
which sends $t_1, t_2$ to either $u(f(t_2) - f(t_1))$ if $t_1 \neq t_2$
or $u(f'(t))$ if $t_1 = t_2 = t$.  We may view $\rho(f)$ as a homotopy
between $ev_1(f)$, which is the restriction to the $t_1 = t_2$ edge, and the
restriction to the $t_1 = 0$ and $t_2 = 1$ edges.  But the restriction to these latter
two edges is canonically 
null-homotopic, since their images lie in the southern hemisphere
of $S^{m-1}$.  Thus, $ev_1$ restricted to $\Emb(\I, \I^m)$ is null-homotopic.
Since $ev_1$ is an equivalence on $\Imm(\I, \I^m)$ this implies that the 
inclusion of $\Emb(\I, \I^m)$ is null-homotopic.

That $\F_m  \simeq \Emb(\I, \I^m) \times \Omega \Imm(\I, \I^m)$ is immediate from 
its definition as the homotopy fiber of this inclusion, and that this is in turn
weakly equivalent to $\Emb(\I, \I^m) \times \Omega^2 S^{m-1}$ follows from the Hirsch-Smale 
theorem \cite{Smal59}.
\end{proof}

The elements of  $\F_3$, which by the above is homotopy equivalent
to $\Emb(\I, \I^3) \times \Omega^2 S^2$, are naturally equipped with framings. 

\begin{proposition}[Proposition 5.14 from \cite{Voli04.1}]\label{P:frame}
The components of $\F_3$ are canonically identified with isotopy classes of framed
knots whose framing is even.
\end{proposition}

Given an element of $\F_3$, a knot with a homotopy through immersions to the unknot, 
Volic ``carries''  the zero framing on the unknot through
the homotopy to define a framing
number on the knot, which is necessarily even.  Note that $\F_m$ is also the
homotopy fiber of the inclusion of the space of framed embeddings in the space
of framed immersions.  When $m=3$ the space of framed immersions is homotopy 
equivalent to $\Omega SO(3)$, so the long-exact sequence in homotopy for this
fibration reads
$$ \cdots \pi_2(SO(3)) = 0 \to \pi_0(\F_3) \to \pi_0 \Emb(\I, \I^3) \times \Z \to 
\pi_1(SO(3)) \cong \Z/2 \to 0,$$
consistent with the calculation that $\pi_0 (F_3) \cong \pi_0( \Emb(\I, \I^3) \times 2\Z$.  

\section{The main result}\label{S:main}

We assemble our work to this point to prove the main result.
As needed for the calculus of functors, extend $\F_m$ to be a functor on the
open sets of $\I$ by sending $U$ to the homotopy fiber of the inclusion
$\Emb(U, \I^m) \to \Imm(U, \I^m)$.  We will recover models for $\F_m$ from
those for embeddings and immersion spaces.

\begin{lemma}\label{L:fib}
If $A$ and $B$ are two functors from $\U(M)$ to spaces with a natural transformation $\tau$
between them, and $F$ is defined by
$F(U) = \hofib(\tau: A(U) \to B(U))$, then $T_k(F) = \hofib(T_k(A) \to T_k(B))$.
\end{lemma}

\begin{proof}
The equality is immediate from the definition of $T_k$, since taking 
homotopy fibers commutes with taking homotopy limits.
\end{proof}

In defining a fiber to  $\rho^m_k$ we are led to the following.

\begin{definition}
Let $e_i : A_n \la \I^m \ra \to A_{n+1} \la \I^m \ra$ send $(u_{j\ell})$ to $(v_{j \ell})$ where
$v_{i, i+1} = *_S$, the basepoint of $S^{m-1}$ and other $v_{j \ell}$ are equal to 
$u_{\sigma_i(j) \sigma_i(\ell)}$.  As before $\sigma_i(j) = j$ or $j - 1$ if $j < i$ or $j > i$
respectively.  By abuse, use $e_i$ to denote its restriction to $C_n \la \I^m, \partial \ra$
mapping to $C_{n+1} \la \I^m, \partial \ra$, as one can check  using \refT{cnlasub2}.
\end{definition}

Alternately, $e_i : C_n \la \I^m, \partial \ra \to C_{n+1} \la \I^m, \partial \ra$ is the 
restriction of $\delta^i$ to $C_n \la \I^m, \partial \ra \times (*_S)^n \subset
C_n' \la \I^m, \partial \ra$.

\begin{definition}
Let $\mF^m_k$ be the functor from $P_0(k)$ to spaces which sends
$S \subseteq [k]$ to $C_{\#S - 1} \la \I^m, \partial \ra$ and which sends
the inclusion $S \subset S \cup j$ to the map $e^i$ where $i$
is the number of elements of $S$ less than $j$.
Let $F^m_k = \holim{}{ \mF^m_k}$.
\end{definition}

\begin{theorem}\label{T:mainpt1}
$F^m_k$ is homotopy equivalent to $T_k \F_m$.  For $m>3$, 
$\eta_\infty : \F_m \to \holim{} {T_k \F_m}$
is a weak equivalence.
\end{theorem}

\begin{proof}
We use the models $D^m_k$ and $G^m_k$
for $T_k \Emb$ and $T_k \Imm$ as given in 
Theorems~\ref{T:modemb}~and~\ref{T:modimm} respectively.
By \refP{emimsq}, $p^m_k : D^m_k \to G^m_k$ agrees with $T_k$ of the inclusion
from embeddings to immersions.  Applying  \refL{fib} 
with $A = \Emb(-, \I^m)$, $B = \Imm(-, \I^m)$, and the natural 
transformation between them be the  standard inclusion, 
we have that $T_k \F_m = \hofib p^m_k$.  

If a map of diagrams  indexed by $P_0(k)$ is a fibration object-wise, 
then the induced map on homotopy limits is a fibration and the fiber is
given by the homotopy limit of the fibers object-wise 
(see for example Lemma~3.5 of \cite{BCSS03}).
Because $\rho^m_k$ is a fibration object-wise, we identify $\hofib p^m_k
= \hofib (\holim{}{\rho^m_k})$ with such a homotopy limit
of object-wise fibers.  By our definition, the diagram of
fibers is $\mF^m_k$, whose homotopy limit is  $F^m_k$, 
establishing the first half of the theorem.

The second half of the theorem is immediate from 
Theorems~\ref{T:imm}~and~\ref{T:GK}.
\end{proof}

Because of \refP{null}, we could alternately
extend $\F_m$ to a functor on $\U(\I)$ by setting $\F^!_m(U) =
\Emb(U, \I^m) \times \Omega \Imm(U, \I^m)$.  The extension $\F^!_m$ would
lead to a set of approximations to $\F_m$  different from the $F^m_k$.

Recall  \refP{assoctoK} that $\K_m$ is an operad with multiplication, which
using \refD{hoch} has an associated cosimplicial object.  We translate from 
$F^m_k$ to $\tot(\K_m^\bullet)$, essentially through the standard projection 
$C_n \la \I^m, \partial \ra \to \wt{C_n} \la \R^m \ra$.  We modify
both $C_n \la \I^m, \partial \ra$ and this projection to define a natural transformation.

\begin{definition}\label{D:witheps}
\begin{itemize}
\item Let  $\varepsilon \leq \frac{1}{6}$. 
For $x \in \R^m$, let $d_+ (x)$ be the distance in $\R^m$ from $x$ to 
$*_+ = (0, \ldots, 0,1)$ and $d_- (x)$ be the distance to $*_-$.
Let $\gamma_j : \R^m \to \R$ be projection onto the $j$th coordinate.
\item  Let $C_n \la \I^m, \partial_\varepsilon \ra$ be the
subspace of $(x_i) \times (u_{ij}) \in C_n \la \I^m, \partial \ra$ where
if $d_+ (x_i)$ and $d_+ (x_j)$ are less than $\varepsilon$ and $i < j$
then $\gamma_k(x_i) = \gamma_k(x_j)$ for $k < m$ and $\gamma_m(x_i) \geq \gamma_m(x_j)$.
Moreover, if $x_i = x_j$ and $i<j$ then $u_{ij} = *_S$.
\item Let $\mF^m_{k, \varepsilon}$ be the functor from $P_0(k)$ to spaces which sends
$S \subseteq [k]$ to $C_{\#S - 1} \la \I^m, \partial_\varepsilon \ra$ and which sends
the inclusion $S \subset S \cup j$ to the map $e^i$ where $i$
is the number of elements of $S$ less than $j$.
\item Let $F^m_{k, \varepsilon} = \holim{}{ \mF^m_{k, \varepsilon}}$.
\end{itemize}
\end{definition}

\begin{proposition}\label{P:step1}
The map $F^m_{k, \varepsilon} \to F^m_k$, induced by the natural transformation
$\iota : \mF^m_{k, \varepsilon} \to \mF^m_{k}$ which at each entry is
the canonical inclusion, is a homotopy equivalence.
\end{proposition}

\begin{proof}
It suffices to check $\iota$ is a homotopy equivalence object-wise, for which
we adapt the machinery developed in \cite{Sinh03} for compactified
configuration spaces.  
Both  $C_{k} \la \I^m, \partial_\varepsilon \ra$ and  $C_k \la \I^m, \partial \ra$
are quotients of the canonical Axelrod-Singer compactifications
which we call $C_k [\I^m, \partial_\varepsilon]$ and $C_k [\I^m, \partial]$
respectively; see Definitions~1.3~and~4.18 of \cite{Sinh03} for the definition
of $C_k [\I^m, \partial]$, which can be modified as in \refD{witheps} for
$C_k [\I^m, \partial_\varepsilon]$.  These quotient maps are
homotopy equivalences, by the proof of Theorem~5.10 of 
\cite{Sinh03}, which applies verbatim in these cases.  

$C_k [\I^m, \partial_\varepsilon]$ retracts 
to its subspace $C_k[\I^m - N^\varepsilon_{\pm}]$,
where $N^\varepsilon_\pm$ is the union of the $\varepsilon$ neighborhoods of $*_+$
and $*_-$ by scaling the $x_i$ by $1-\varepsilon$.
Both $C_k [\I^m, \partial]$ and $C_k[\I^m - N^\varepsilon_{\pm}]$
are manifolds with corners (see Theorem~4.4 of \cite{Sinh03}), 
and thus are homotopy equivalent
to their interiors, $C_k(Int(\I^m))$ and $C_k(Int(\I^m - N^\varepsilon_{\pm}))$
respectively.   But these interior configuration spaces
are diffeomorphic, since $Int(\I^m)$ and $Int(\I^m - N^\varepsilon_{\pm})$
are.  Composing this diffeomorphism with the previous homotopy
equivalences establishes the equivalence of 
$C_k \la \I^m, \partial_\varepsilon \ra$ and $C_k \la \I^m, \partial \ra$ and
thus establishes the result.

\end{proof}

We use $C_k \la \I^m, \partial_\varepsilon \ra$ because
they readily project to $\wt{C_k}\la \R^m \ra$ in a way compatible with
its cosimplicial structure maps.

\begin{definition}\label{D:pk}
\begin{itemize}
\item Let $(a_i)_{i=1}^m$, $a_i \in \R$  denote a point in $\R^m$. 
Define $\lambda_+ : (\R^m - *_+) \to \R^m$ by 
sending $(a_i)$ to $(b_i)$ where if $i \neq m$ then $b_i = a_i$
and $$b_m = 
\begin{cases}
\frac{\varepsilon a_m}{d_+(a_i)} \;\; d_+(a_i) < \varepsilon \\
a_m             \;\;\;\;\;\;  d_+(a_i) \geq \varepsilon.
\end{cases}
$$
Define $\lambda_- : (\R^m -  *_-) \to \R^m$ similarly, and let 
$\lambda = \lambda_+ \circ \lambda_-$.
\item Define $\pi_k : C_k \la \I^m, \partial_\varepsilon \ra \to \wt{C_k} \la \R^m \ra 
\subset (S^{m-1})^{\binom{k}{2}}$ 
by sending $(x_i) \times (u_{ij})$ to $(v_{ij})$
where $v_{ij}$ is:
\begin{itemize}
\item $u(\lambda(x_i) - \lambda(x_j))$ if $x_i \neq x_j$ 
and neither equals $*_+$ or $*_-$.
\item The Jacobian on $\lambda$ applied to $u_{ij}$ if $x_i = x_j$.
\item $*_S$, if $x_i \neq x_j$ and either $x_i = *_+$ or $x_j = *_-$.
\end{itemize}
\end{itemize}
\end{definition}

\begin{proposition}
$\pi_k$ is continuous.
\end{proposition}

\begin{proof}
We first identify $\pi_k$ on the subspace t $C_k\la \I^m - (*_+ \cup *_-)\ra$ with the
composite of $C_k \la \lambda \ra$, the map on configuration spaces
induced by the embedding $\lambda$ (see Corollary~4.8 of \cite{Sinh03}), and 
the canonical projection
$C_k \la \R^m \ra$ to $\wt{C_k} \la \R^m \ra$.  What remains is to 
check continuity on the subspace
in which some $x_n = *_+$.  Consider a sequence 
$\{(x^\ell_i), (u^\ell_{ij})\}_{\ell = 1}^\infty$ with limit point $(x_i^\infty) \times (u^\infty_{ij})$,
so that $x_n^\infty = *_+$.  We show that its image under $\pi_k$ has $v_{nj}$
which approaches $*_S$ if $x_j^\infty \neq *_+$ or $n < j$ or which approaches $-*_S$ otherwise.  
For each $j$, either $x^\infty_j \in N^\varepsilon_+$, 
which is also true for $\ell$ sufficiently large, in which case
$u_{nj}$ must be $*_S$ if $n<j$ or  $-*_S$ if $n>j$, so 
that the sequence $v_{nj}$ would be eventually constant at $*_S$
or $-*_S$. Or  if $x^\infty_j \notin N^\varepsilon$ then as $x_n^\ell \mapsto *_+$ the last coordinate
of $\lambda_+({x_i}^\ell)$ becomes arbitrarily large.  Because  ${x_j}^\ell \mapsto x_j$ stays 
in $\I^m$ we have $u(\lambda(x^\ell_n) - \lambda(x^\ell_j)) \mapsto  *_S$.
Continuity when some $x_n = *_-$ works similarly.
\end{proof}

We now may assemble our main result, \refT{main}, which casts the embedding
calculus tower for $\F_m$ in the language of operads.   For convenience, we restate
the theorem here.

\begin{theorem}\label{T:main2}
The $k$th approximation to $\F_m$ in the embedding calculus, namely
$T_k \F_m$, is weakly equivalent to $\tot^k \K_m^\bullet$.  
\end{theorem}

\begin{proof}
We will check that the maps $\pi_k$ assemble to a natural transformation of 
functors from $\mF^m_k$ to $ \K_m^\bullet \circ c_k$, with $c_k$ as in 
\refD{ck}, which gives rise to
a weak equivalence on homotopy limits.  
\refT{mainpt1} then says that the homotopy limit of $\mF^m_k$ is
weakly equivalent to $T_k \F_m$, and \refL{coshol}
implies that the homotopy limit of $\K_m^\bullet \circ c_k$ is weakly
equivalent to $\tot^k \K_m^\bullet$, establishing the theorem.

For the assembled $\pi_k$ to be a natural transformation, 
we must have $\pi_k \circ e^i = d^i \circ \pi_k$.  For most $i$ this is 
immediate to check, as repeating coordinates and passing to the
quotient $\wt{C}_k \la \R^m \ra$ are processes which clearly commute.
The $i=0$ and $i=k+1$ cases require the modifications we made in \refD{pk}.
For $\K_m^\bullet \circ c_k$ we trace through 
Definitions~\ref{D:hoch}~and~\ref{D:mB}, \refT{kop} and \refP{assoctoK}
to see that $d^{k+1}$ takes a point $(u_{ij}) \in \wt{C}_k \la \R^m \ra$, 
leaves all these $u_{ij}$ unchanged, and adds $u_{i, k+1} = *_S$
for all $i$ to obtain a point in $ \wt{C}_{k+1} \la \R^m \ra$.  On
the other hand, $e^{k+1}$ adds the $k+1$st point to the configuration at 
$*_-$, which under $\pi_k$ will also lead to all $u_{i, k+1} = *_S$.  The
$i=0$ case works similarly.

The fact that the assembled $\pi_k$ induce a weak equivalence
on homotopy limits follows from it being a homotopy equivalence object-wise.
We already know from the proof of \refP{step1} that 
$C_k \la \I^m, \partial_\varepsilon \ra$ is homotopy equivalent to the
subspace $C_k (Int(\I^m - N^\varepsilon_{\pm}))$, which is diffeomorphic
to $C_k(\R^m)$.  Composed with this diffeomorphism 
on this subspace, $\pi_k$ is the standard projection
$C_k(\R^m) \to \wt{C_k}(\R^m)$ followed by the canonical map to 
$\wt{C_k} \la \R^m \ra$ which is a homotopy equivalence by 
Corollaries~4.5~and~5.9 of \cite{Sinh03}.
\end{proof}

The first half of \refC{main}, which states that for $m>3$ the totalization of the
Kontsevich operad faithfully models the weak homotopy type of $\F_m$,
now follows from \refT{mainpt1}.  The second half of \refC{main}, which states that
for $m=3$ all real-valued finite-type framed knot invariants are pulled back from this 
operad model, follows from the
main results of \cite{Voli04.1}.  We may take any framed knot, double its framing, and
then apply \refP{frame} to get a corresponding component of  $\F_3$.  
In \cite{Voli04.1}, Volic uses Bott-Taubes integrals to
define finite-type invariants on the Taylor tower for $\F_3$, which we have now
identified with the Tot tower for $\K_3$.  See in 
particular Theorems~1.2 and 1.3 and Section~6.4 of \cite{Voli04.1}.

\section{Observations and consequences}\label{S:cons}

\subsection{Spectral sequences}

The results in this section parallel those of section 7 of \cite{Sinh02}.
Applying the homotopy  spectral sequence of \refP{csshot} 
for $\K_m^\bullet$ we immediately have the following.

\begin{theorem}\label{T:knothot}
There is a spectral sequence converging to $\pi_*(\tot \K_m^\bullet)$
with $$E_1^{-p,q} = \bigcap \; {\text{ker}} \; {s^k}_* \subseteq 
\pi_q(C_p (\R^m)). $$ The $d_1$ differential is the restriction to 
this kernel of the map $$\Sigma_{i=0}^{p+1} (-1)^i {d^i}_* \colon 
\pi_q(C_{p-1} (\R^m)) \to \pi_q(C_p(\R^m)).$$ 
\end{theorem}

\refT{mainpt1} implies that this spectral 
sequence computes homotopy groups of $\F_m$ when $m \geq 4$.
Except for in the
$p = 1$ column, this spectral sequence coincides exactly with that studied
with rational coefficients in \cite{ScSi02}, so we do not give a
more explicit description here.  Kontsevich \cite{Kont00} has also 
examined the rows of this spectral sequence.

For $m=3$, the case of classical knots, we conjecture that 
$\eta_k : \F_m \to T_k \F_m$ is a universal type-$(k-1)$ framed
knot  invariant over the integers.  For $k \leq 3$, we may deduce 
this from the main results of  \cite{BCSS03}.   
In unpublished work, Conant has shown that the entries $E^2_{-k,k}$ of this
spectral sequence are isomorphic to the module of primitive
weight systems of degree $k-1$ over the integers,  a
first step to this conjecture in full generality.

In light of \refT{poisson}, the homology spectral sequence from
\refT{csshom} has a pleasant description.   
Recall \refD{hoch}, which for operads of vector spaces introduces
the notation of $HH^*(\Op)$ for the total cohomology of the 
associated cosimplicial object.

\begin{theorem}\label{T:HHss}
There is a spectral sequence with $E^2_{-p,q} = HH^{p,q}({\rm Poiss}_m)$
which for $m \geq 4$ converges to the homology of $\tot \K^\bullet_m$, and
thus of $\F_m$.
\end{theorem}

\begin{proof}
If we use the second description of the homology spectral sequence
from \refT{csshom}, then $E^1_{-*,*}$ will be $H_*(\K_m^\bullet)$, 
which is the Poisson operad by \refT{poisson}.  The induced
operad with multiplication structure on the Poisson operad
is the standard one.  Thus, the $d^1$ differential will coincide with
the differential for total cohomology of the Poisson operad, and the
$E^2$ term will be the total (or Hochschild) cohomology
of the Poisson operad as stated.

It remains to check the convergence conditions of \refT{csshom}.  
In the case of the Kontsevich operad, the entries $\K_m^k = \wt{C}_k \la \R^m \ra$ are
homotopy equivalent to $C_k(\R^m)$, which are
simply connected if $m \geq 3$.
Using the first definition of \refT{csshom}, we start with $H_*(C_p(\R^m))$
and explicitly understand the kernels of the maps $s^i_*$.  We use
\refT{poisson} and \refD{poiss} to identify $H_*(C_p(\R^m))$ in terms of products
of brackets in variables $x_1, \ldots, x_p$.  Tracing through the definitions
of the associated cosimplicial object, $s^i$ sends a product of brackets
in the $x_j$ to either zero, if the variable $x_i$ appears in a bracket, or
the monomial obtained by removing $x_i$ and relabeling $x_j$ to $x_{j-1}$ for $j>1$,
if $x_i$ does not appear in a bracket.  Therefore to be in the kernel of all of the $s^i$,
all of the variables $x_i$ must appear in a bracket, so there must be at least $\frac{k}{2}$
brackets, leading to a total degree of at least $\frac{k(m-1)}{2}$.  For $m > 3$, this 
is greater than $k$ and thus gives the estimate needed for application of \refT{csshom}.
\end{proof}

This spectral sequence in real-valued cohomology coincides with the
homotopy spectral sequence of the Taylor tower for the functor to spectra
$U \mapsto \R \wedge E_m(U)$, the real Eilenberg-MaClane spectrum
smashed with $E_m(U)$.  For $m=3$, Volic's results \cite{Voli04.1, Voli04.2}
imply that the map from the knot space to this Taylor tower serves
as a universal framed finite-type invariant over the real numbers.

\subsection{A little two-cubes action from the McClure-Smith framework}

\refT{main} fits perfectly into the framework created by 
McClure and Smith in their
solution of the Deligne conjecture over the integers \cite{McSm01}.  
One of their central results is the following.

\begin{theorem}\label{T:mcsm}
The totalization of the associated cosimplicial object of an operad with
multiplication
admits an action of an operad equivalent to the little 2-cubes operad,
as does its homotopy invariant totalization.
\end{theorem}

\begin{proof}
We are simply collecting results from \cite{McSm01} and \cite{McSm02}.
For the standard totalization, we are simply quoting Theorem 3.3 in \cite{McSm01}.
For the homotopy invariant totalization, Theorem~15.3 of \cite{McSm02} says that 
$\tot$ of any cosimplicial space with what is called a $\Xi^2$-structure has
an action of an operad equivalent to the little 2-cubes.  
Proposition~10.3 of \cite{McSm02} identifies an operad with multiplication structure
on a sequence of spaces with a $\Xi^2$ structure.  
\end{proof}

\begin{example}\label{X:omega2}
Consider the cosimplicial model for the space of maps from $S^2$ to $X$,
namely $X^{S^2_\bullet}$.  By \refT{binoms2} $S^2_\bullet \cong \mB^\bullet$, 
so there is an operad structure on this collection of spaces.  In order to get
an operad with multiplication, we restrict each $X^{\mB^n} = (x_\alpha)$
to the subspace in which $x_+ = *$, where $+$ is the basepoint of $\mB^n$
and $*$ is the base point of $X$.  The operad structure maps restrict appropriately,
and we obtain ${X}_\star^{S^2_\bullet}$, to which the associative operad maps
at each level to the point with all $x_\alpha = *$.

Applying \refT{mcsm}, the totalization of ${X}_\star^{S^2_\bullet}$ is a 
little 2-cubes space, and we know that its totalization is $\Omega^2 X$.  
McClure and Smith 
fully develop this example (more generally for the standard model of $\Omega^n X$) 
in Section~11  of \cite{McSm02}.  They show that the little 2- (or $n-$) cubes action which arises
in these examples coincide
with the standard ones. 
\end{example}


We can immediately establish \refT{2cube}, one of our main results, 
which we restate here.

\begin{theorem}
For any $m$, there is a little two-cubes action on $\tot(\K^\bullet_m)$.  For $m>3$,
$\F_m$ is a two-fold loop space.
\end{theorem}

\begin{proof}
Applying \refT{mcsm} for the Kontsevich operad with its given multiplication 
establishes the two-cubes action.

By \refT{main}, if $m \geq 4$,
$\F_m$ is homotopy equivalent to $\tot(\K_m^\bullet)$,
so it has a 2-cubes action as well.  But $\F_m$ 
is connected for $m \geq 4$, since it is the product of  $\Omega(\Imm(\I, \I^m)) 
\simeq \Omega^2 S^{m-1}$ and $\Emb(\I, \I^m)$ which are both connected
(that the latter space is connected is because any path through
maps from an embedding to the standard one becomes
an isotopy once put in general position).
By the recognition theorem of \cite{Boar71, BoVo73}, 
$\F_m$ is a 2-fold loop space.
\end{proof}

Our operad model for $\F_m$
has  already been important for closer examination
of the homotopy type  of $\F_m$.  
In \cite{LaVo05}, Volic and Lambrechts establishes a formality
result for this model, 
which determines the rational homotopy type of $\F_m$ for $m>3$.

A two-cubes action has been used with spectacular success in dimension
three.  In \cite{Budn03}, Budney constructs a little two-cubes action directly on
the space of framed (long) knots in any dimension.  
He goes on to show that the two-cubes action is free
when $m=3$, generated by the components of prime knots. 
In  \cite{Budn05.1}, he determines the homotopy types of these prime 
components in terms of the JSJ decompositions of prime knots.

In further work, we may develop a
model for framed knots closely related to those studied in 
this paper.  Namely, there is an operad whose entries are
$\wt{C}_n \la \R^m \ra \times (SO(m))^n$, where in the operad composition,
the elements of $SO(m)$ act on the configurations (compare
with the operad of ``turning balls'').  The techniques in this paper 
should also show that the totalization of this operad models spaces
of framed knots.  It may then be fruitful to compare Budney's geometric
two-cubes action to the one which arises from the McClure-Smith machinery,
as well as to combinatorial product and bracket structures on homology
defined by Tourtchine \cite{Tour03}.


\begin{thebibliography}{10} 

\bibitem{AxSi94} S. Axelrod and I. Singer. Chern-Simons perturbation theory, II.  
{\em Jour. Diff. Geom.} {\bf 39} (1994), no. 1, 173--213.   

\bibitem{BCWW04} A. J. Berrick, F. Cohen, Y. Wong and J. Wu.
Configurations, braids, and homotopy groups, preprint 2004.

\bibitem{Boar71} J. M. Boardman.  Homotopy structures and the language
of trees.  In {\em Algebraic Topology}, Proc. Sympos. Pure Math, Vol. XXII,
AMS, 1971.  pp 37--58.

\bibitem{BoVo73} J. M. Boardman and R. Vogt. Homotopy invariant algebraic structures
on topological spaces. Lecture Notes in Mathematics, Vol. 347, 1973.

\bibitem{BoKa72} A. Bousfield and D. Kan.  Homotopy limits, completions 
and localizations.  Lecture Notes in Mathematics, Vol. 304, 1972.

\bibitem{Bous86} A. Bousfield.  On the homology spectral sequence of a
cosimplicial space.  {\em Amer. J. Math.} {\bf 109} (1987), no. 2, 361--394.

\bibitem{Budn03} R. Budney.  Little cubes and long knots.  To appear in 
{\em Topology}.

\bibitem{Budn05.1} \bysame Topology of spaces of knots in dimension 3.  
math.GT/0506524.

\bibitem{BCSS03} R. Budney, J. Conant, K. Scannell and D. Sinha. New perspectives
on self-linking.  {\em Advances in Mathematics}. {\bf 191} (2005) 78--113.

\bibitem{FCoh76} F. Cohen. The homology of $C_{n+1}$ spaces.  In Lecture 
Notes in Mathematics {\bf 533} (1976).

\bibitem{FuMa93} W. Fulton and R. MacPherson.  Compactification  of configuration spaces.   
{\em Annals of Mathematics} {\bf  139} (1994), 183--225. 

\bibitem{Gaif03} G. Gaiffi.  Models for real subspace arrangements and stratified
manifolds.   {\em Int. Math. Res. Not.}  (2003)  no. 12, 627--656.

\bibitem{GeVo94} M. Gerstenhaber and A. Voronov,  Homotopy $G$-algebras and
moduli space operads.  {\em Intern. Math. Res. Notices} {\bf 3} (1995) 141--153.

\bibitem{GeJo94} E. Getzler and J. Jones.  Operads, homotopy algebra and iterated integrals 
for double loop spaces.  hep-th/9403055.

\bibitem{Good91}T. Goodwillie. Calculus II, Analytic functors. 
$K$-Theory {\bf 5} (1991/92), no. 4, 295--332.

\bibitem{GoWe99} T. Goodwillie and M. Weiss.  Embeddings from the point of 
view of immersion theory II. {\em Geometry and Topology} 3 (1999), 
103--118. 

\bibitem{GKW01} T. Goodwillie, J. Klein and M. Weiss.   
Spaces of smooth embeddings, disjunction and surgery. Surveys on 
surgery theory, Vol. 2, 221--284, Ann. of Math. Stud., {\bf 149} (2001). 

\bibitem{GKW02} T. Goodwillie, J. Klein and M. Weiss. 
A Haefliger type description of the embedding calculus tower. 
{\em Topology} 42 (2003), no. 3, 509--524.

\bibitem{GoKl03} T. Goodwillie and J. Klein.  \newblock Excision 
statements for spaces of smooth embeddings.  \newblock In preparation.

\bibitem{Good03} T. Goodwillie.  Excision estimates for spaces of homotopy
equivalences.  Preprint is available at: {\tt http://www.math.brown.edu/faculty/goodwillie.html}.

\bibitem{Hirs03}
P.~S. Hirschhorn, \emph{Model categories and their localizations}, Math.
  Surveys Monographs, vol.~99, Amer. Math. Soc., Providence, 2003.


\bibitem{Kont99} M. Kontsevich, Operads and motives in deformation
quantization, {\em{Lett. Math. Phys.}} 48 (1999) 35--72.

\bibitem{Kont00} \bysame. Operads of little discs in algebra 
and topology, Lecture at the Mathematical Challenges Conference, UCLA, 
2000. 

\bibitem{LaVo05} P. Lambrechts and I. Volic.  The rational homotopy type of
spaces of knots in high codimension.  In preparation.

\bibitem{Mark99} M. Markl.   A compactification of the real configuration space as an 
operadic  completion.  {\em  J. Algebra}  {\bf 215}  (1999),  no. 1, 185--204.

\bibitem{MSS02} M. Markl, S. Shnider, and J. Stasheff. 
{\em Operads in algebra, topology and physics}. Math. Surv. and Monographs, 96. 
AMS, Providence, RI, 2002.

\bibitem{May72} J. P. May. The geometry of iterated loop spaces. Lecture Notes in 
Mathematics, Vol. 271, 1972.

\bibitem{McSm01}  J. McClure and J. Smith. 
 A solution of Deligne's Hochschild cohomology conjecture.  In:
Recent progress in homotopy theory (Baltimore, MD, 2000), 153--193,
Contemp. Math., {\bf 293}, AMS, 2002. 

\bibitem{McSm02}
\bysame  Cosimplicial objects and little $n$-cubes. I.  math.QA/0211368.

\bibitem{McSm04}
\bysame 
Operads and cosimplicial objects: an introduction.
In {\em Axiomatic, Enriched and Motivic Homotopy Theory.}
Edited by J.P.C. Greenlees. Kluwer 2004.

\bibitem{Misl95} G. Mislin. Wall's finiteness obstruction.  In {\em Handbook of algebraic
topology}, edited by I. M. James.  Elsevier, 1995.


\bibitem{Pala66}
R. Palais.  Homotopy theory of infinite-dimensional manifolds. 
{\em Topology}  {\bf 5} (1966) 1-16.

\bibitem{Rect70} D. Rector.  
Steenrod operations in the Eilenberg-Moore spectral sequence. 
{\em Comment. Math. Helv.} {\bf 45} 1970 540--552. 

\bibitem{ScSi02} K. Scannell and D. Sinha. A one-dimensional embedding complex.
{\em Journal of Pure and Applied Algebra}  {\bf 170} (2002) No. 1, 93--107. 

\bibitem{Ship96} B. Shipley. Convergence of the homology spectral sequence  
of a cosimplicial space. Amer. J. Math. {\bf 118} (1996), no. 1, 179--207. 

\bibitem{Sinh02} D. Sinha.  The topology of spaces of knots. math.AT/0202287.

\bibitem{Sinh03} \bysame Manifold theoretic compactifications of
configuration spaces.  To appear in {\em Selecta Math.}.  math.GT/0306385.

\bibitem{Sinh04.2} \bysame Configuration spaces, Hopf invariants, and Whitehead products.
In preparation.

\bibitem{Smal59} S. Smale. The classification of immersions of spheres 
in Euclidean spaces.   {\em Annals of Mathematics} (2) {\bf 69} 1959 327--344. 

\bibitem{Tour03} V. Tourtchine.  On the homology of long knots.  To appear
in NATO Conference Proceedings. math.QA/0105140.

\bibitem{Vass92}  V. Vassiliev. 
\newblock {\em Complements of discriminants of smooth maps:  
topology and applications.} 
Translations of Mathematical Monographs, 98. American 
Mathematical Society, Providence, RI, 1992. 

\bibitem{Vass97} \bysame  Topology of two-connected graphs and homology of spaces of knots.  In: Differential and symplectic topology of knots and curves,  253--286, AMS Transl. Ser. 2, {\bf 190},  1999.


\bibitem{Voli04.1} I. Volic.  Finite type knot invariants and calculus of functors. 
To appear in {\em Compositio Mathematica}.  math.AT/0401440.

\bibitem{Voli04.2} \bysame  Configuration space integrals and Taylor towers for spaces of knots. 
math.GT/0401282.

\bibitem{Weis95} M. Weiss.  Calculus of embeddings.  
{\em Bulletin of the AMS} 33  (1996),  177-187. 


\bibitem{Weis99} 
\bysame Embeddings from the point of view of immersion theory I.  
{\em Geometry and Topology} 3 (1999), 67--101. 

\end{thebibliography}
\end{document}